\documentclass[12pt]{amsart}
\usepackage{amsmath,amscd,amssymb,amsfonts,graphics,mathrsfs}
\usepackage{xcolor,hyperref}
\hypersetup{colorlinks=true,citecolor=black,linkcolor=black}
\newcommand{\hl}{\hyperlink}
\newcommand{\htt}{\hypertarget}
\newcommand{\h}{\hbox}
\newcommand{\q}{\quad}
\newcommand{\nin}{\noindent}
\newcommand{\bs}{\par\bigskip}
\newcommand{\ms}{\par\medskip}
\newcommand{\sk}{\par\smallskip}
\newcommand{\bsn}{\par\bigskip\noindent}
\newcommand{\msn}{\par\medskip\noindent}

\newcommand{\ges}{\geqslant}
\newcommand{\les}{\leqslant}
\newcommand{\1}{\hskip1pt}

\newcommand{\mopl}{\hbox{$\bigoplus$}}
\newcommand{\msum}{\hbox{$\sum$}}
\newcommand{\mprod}{\hbox{$\prod$}}
\newcommand{\B}{{\mathscr B}}
\newcommand{\Cc}{{\mathscr C}}
\newcommand{\D}{{\mathscr D}}
\newcommand{\E}{{\mathscr E}}
\newcommand{\G}{{\mathscr G}}
\newcommand{\Hc}{{\mathscr H}}

\newcommand{\M}{{\mathscr M}}

\newcommand{\OO}{{\mathscr O}}
\newcommand{\Mf}{{\mathscr M}_{\!_f}}
\newcommand{\Mtf}{{}\,\,\,\,\,\widetilde{\!\!\!\!\!\mathscr M}{}_{\!\!_f}}
\newcommand{\Nf}{{\mathscr N}\!\!_f}
\newcommand{\bt}{\widetilde{b}}
\newcommand{\ct}{\widetilde{c}}
\newcommand{\Ft}{\widetilde{F}}

\newcommand{\Rt}{\widetilde{R}}
\newcommand{\Zt}{\widetilde{Z}}
\newcommand{\delt}{\widetilde{\delta}}
\newcommand{\alt}{\widetilde{\alpha}}
\newcommand{\nut}{\widetilde{\nu}}
\newcommand{\Yt}{\widetilde{Y}}
\newcommand{\ft}{\widetilde{f}}
\newcommand{\gt}{\widetilde{g}}
\newcommand{\RR}{{\mathbf R}}

\newcommand{\PP}{{\mathbb P}}
\newcommand{\Q}{{\mathbb Q}}
\newcommand{\C}{{\mathbb C}}
\newcommand{\N}{{\mathbb N}}
\newcommand{\DD}{{\mathbf D}}
\newcommand{\Z}{{\mathbb Z}}

\newcommand{\Hppf}{\mathscr{H}''_{\!f}}
\newcommand{\Htppf}{\widetilde{\mathscr{H}}''_{\!f}}
\newcommand{\rf}{r\!_f}
\newcommand{\PD}{P^{\D}}
\newcommand{\Gr}{{\rm Gr}}
\newcommand{\DR}{{\rm DR}}
\newcommand{\la}{\lambda}
\newcommand{\al}{\alpha}
\newcommand{\be}{\beta}
\newcommand{\om}{\omega}
\newcommand{\Om}{\Omega}
\newcommand{\lng}{\langle}
\newcommand{\rng}{\rangle}
\newcommand{\dd}{\partial}
\newcommand{\ddd}{{\rm d}}
\newcommand{\eq}{\,{=}\,}
\newcommand{\defs}{\,{:=}\,}
\newcommand{\gess}{\,{\ges}\,}
\newcommand{\less}{\,{\les}\,}
\newcommand{\sgt}{\,{>}\,}
\newcommand{\slt}{\,{<}\,}
\newcommand{\nes}{\,{\ne}\,}
\newcommand{\mi}{\1{-}\1}

\newcommand{\pl}{\1{+}\1}
\newcommand{\scd }{\1{\cdot}\1}
\newcommand{\bl}{\bigl}
\newcommand{\br}{\bigr}
\newcommand{\sst}{\,{\subset}\,}
\newcommand{\stm}{\,{\setminus}\,}
\newcommand{\ins}{\,{\in}\,}
\newcommand{\tos}{\,{\to}\,}
\newcommand{\cols}{\,{:}\,}
\newcommand{\ssc}{\,\raise.15ex\hbox{${\scriptstyle\circ}$}\,}
\newcommand{\ssb}{\raise.15ex\h{${\scriptscriptstyle\bullet}$}}
\newcommand{\into}{\hookrightarrow}
\newcommand{\onto}{\twoheadrightarrow}
\newcommand{\simto}{\,\,\rlap{\hskip1.3mm\raise1.4mm\hbox{$\sim$}}\hbox{$\longrightarrow$}\,\,}
\newcommand{\indlim}{\rlap{\raise-6.2pt\h{$\,\to$}}{\rm lim}}
\begin{document}
\h{}\bs
\centerline{\large Length of ${\mathscr D}_Xf^{-\alpha}$ in the Isolated Singularity Case}
\bs
\centerline{Morihiko Saito}
\bsn\ms
\vbox{\nin\narrower\smaller
{\bf Abstract.} Let $f$ be a convergent power series of $n$ variables having an isolated singularity at 0. For a rational number $\alpha$, setting $(X,0)\eq({\mathbb C}^n,0)$, we show that the length of the ${\mathcal D}_X$-module ${\mathcal D}_Xf^{-\alpha}$ is given by $\widetilde{\nu}_{\alpha}\pl r_{\!f}\widetilde{\delta}_{\alpha}\pl 1$. Here $r_{\!f}$ is the number of local irreducible components of $f^{-1}(0)$ (with $r_{\!f}\eq1$ for $n\sgt 2$), $\widetilde{\nu}_{\alpha}$ is the dimension of the graded piece ${\rm Gr}_V^{\alpha}$ of the $V$-filtration on the saturation of the Brieskorn lattice modulo the image of $N:=\partial_tt\mi\alpha$ on ${\rm Gr}_V^{\alpha}$ of the Gauss-Manin system, and $\widetilde{\delta}_{\alpha}\defs 1$ if $\alpha\ins{\mathbb Z}_{>0}$, and 0 otherwise. This theorem can be proved also by employing a generalization a recent formula of T.\,Bitoun in the integral exponent case. The theorem generalizes an assertion by T.\,Bitoun and T.\,Schedler in the weighted homogeneous case where the saturation coincides with the Brieskorn lattice and $N\eq0$. In the semi-weighted-homogeneous case, our theorem implies some sufficient conditions for their conjecture about the length of ${\mathcal D}_Xf^{-1}$ to hold or to fail.}
\ms\bs
\centerline{\bf Introduction}
\bsn
Let $f\ins\C\{x\}$ be a convergent power series of $n$ variables having an isolated singularity at 0 with $f(0)\eq0$. To calculate the monodromy on the vanishing cohomology in an algebraic way, E.\,Brieskorn \cite{Br} introduced the {\it Brieskorn lattice\1} $\Hppf$, which is a free $\C\{t\}$-module of rank $\mu_f$, and is endowed with a regular singular {\it Gauss-Manin connection,} where $\mu_f$ is the Milnor number of $f$. B.\,Malgrange \cite{Ma} proved that the local {\it reduced Bernstein-Sato polynomial\1} $\bt_f(s)\defs b_f(s)/(s{+}1)$ coincides with the {\it minimal polynomial\1} of the action of $-\dd_tt$ on the quotient $\Htppf/t\Htppf$ with $\Htppf$ the {\it saturation\1} of $\Hppf$ defined by
\htt{1}{}
$$\Htppf:=\msum_{j\in\N}\,(\dd_tt)^j\Hppf\,\subset\,\G_f,
\leqno(1)$$
\par\nin with $\G_f\defs\Hppf[\dd_t]$ the {\it Gauss-Manin system.} The latter is the {\it localization\1} of $\Hppf$ by the action of $\dd_t^{-1}$, which is well defined on $\Hppf$, see (\hl{1.1.2}{1.1.2}) below. In this paper Bernstein-Sato polynomials are called {\it BS polynomials\1} for short.
\sk
Recently T.\,Bitoun and T.\,Schedler investigated the length of the quotient $\D_X$-modules $\D_Xf^{-\al}/\D_Xf^{-\al+1}$ for $\al\ins\Q$ in the {\it weighted homogeneous\1} case, where $(X,0)\eq(\C^n,0)$. They found that it coincides with the sum of $\delta_{\al,1}$ (Kronecker delta) and the {\it multiplicity of spectral number\1} $n_{f,\al}$, assuming $n\gess 3$ (and using \hl{3.6}{3.6} below), see \cite[Cor.\,1.20 and Thm.\,2.1]{BiSc}. (This assertion follows also from \cite[Cor.\,1]{rp}.) They conjectured that the assertion for $\al\eq1$ would hold in general, see \cite[Conj.\,1.7]{BiSc}. For $\al\nes1$, this fails because of {\it shifts\1} of roots of $b_f(s)$ under $\mu$-constant deformations, although the multiplicity of spectral number $n_{f,\al}$ remains invariant, see \cite{Kat1}, \cite{Kat2}, \cite{Ca}, \cite{sem} and also Remark~\hl{R1.3a}{1.3a} and \hl{R1.2b}{1.2b} below. (It is easy to see that $\D_Xf^{-\al+1}\eq\D_Xf^{-\al}$ if $b_f(-\al)\nes0$, but its converse does not necessarily hold, see \cite{rp}.)
\sk
In general, the multiplicities $n_{f,\al}$ can be defined by
\htt{2}{}
$$n_{f,\al}\eq\dim_{\C}\Gr_V^{\al}(\Hppf\!/t\1\Hppf)\q(\al\ins\Q),
\leqno(2)$$
\par\nin following \cite{Va1}, where $V$ is the $V$-filtration of Kashiwara and Malgrange {\it indexed by\1} $\Q$ (see also \cite{sup}, \cite[\S3.4]{fil}) so that the action of $N:=\dd_tt{-}\al$ is {\it nilpotent\1} on $\Gr_V^{\al}$, see \hl{1.1}{1.1} and Remark~\hl{R1.2a}{1.2a} below. (This filtration was originally indexed by $\Z$, see \cite{sup}, \cite[\S3.4]{fil} about the reason for which it must be indexed by $\Q$.)
\sk
In the weighted homogeneous case, the $V$-filtration can be induced by the {\it weighted degree\1} of differential forms and we have the equality $\Htppf\eq\Hppf$ (see Remarks~\hl{R1.1}{1.1} and \hl{R1.2c}{1.2c} below). Hence the spectral numbers of $f$ coincide with the roots of the reduced BS polynomial up to sign forgetting the multiplicities.
\sk
In general, set $\Htppf{}^{(\al)}:=\Gr_V^{\al}\Htppf$, $\,\G_f^{(\al)}:=\Gr_V^{\al}\G_f$, and
\htt{3}{}
$$\nut_{f,\al}:=\dim_{\C}{\rm Im}\bl(\Htppf{}^{(\al)}\to\G_f^{(\al)}\!/N\G_f^{(\al)}\br),
\leqno(3)$$
\par\nin see \hl{3.7}{3.7} below for the division by the image of $N$. If $f$ is {\it weighted homogeneous,} we have
\htt{4}{}
$$\nut_{f,\al}=\msum_{j\ges 0}\,n_{f,\al-j},
\leqno(4)$$
\par\nin since $\Htppf\eq\Hppf$ and $N\eq0$.
\sk
Let $\rf$ be the number of local irreducible components of $Z\defs f^{-1}(0)\sst X$. Note that $\rf\eq1$ if $n\gess 3$. Set $\delt_{\al}\defs 1$ if $\al\ins\Z_{>0}$, and $0$ otherwise. We denote by $\ell_{\D_{X,0}}(\M)$ the length of a holonomic $\D_{X,0}$-module $\M$ in general. We prove in this paper the following.
\par\htt{T1}{}\msn
{\bf Theorem~1.} {\it For $\al\ins\Q$, there is the equality
\htt{5}{}
$$\ell_{\D_{X,0}}(\D_{X,0}f^{-\al})=\nut_{f,\al}\pl\rf\delt_{\al}\pl 1,
\leqno(5)$$
\par\nin where $\nut_{f,\al}$, $\rf\delt_{\al}$, $1$ are the lengths of regular holonomic $\D_{X,0}$-modules supported respectively on $0$, local irreducible components of $(Z,0)$, and $(X,0)$.}
\ms
This follows easily from \cite[Thm.\,1 and Prop.\,1]{rp}, that is, Theorem~\hl{T1.4}{1.4}, Proposition~\hl{P1.4}{1.4} in this paper, see \hl{2.1}{2.1} below. There is another approach generalizing the last assertion in \cite[Cor.\,3.0.1]{Bit} to the non-integral exponent case (see Corollary~\hl{C4.1}{4.1} below), which implies a simpler proof of Theorem~\hl{T1}{1}, see \hl{4.2}{4.2} below.
\sk
Let $f_1\ins\C[x]$ be a weighted homogeneous polynomial with an isolated singularity at $0$. Let $f_{>1}\eq\msum_{\be>1}\,f_{\be}$ with $f_{\be}$ weighted homogeneous polynomials of weighted degrees $\be$ for any $\be\gess 1$. Set $f\eq f_1\pl f_{>1}$. This is called a {\it semi-weighted-homogeneous\1} deformation of $f_1$. It is essentially equivalent to that $f$ is a $\mu$-constant deformation of $f_1$ (adding coefficients of monomials as parameters), see \cite{Va3}, where we may assume that monomials form a basis of $V^{>\alt_f+1}(\C\{x\}/(\dd f_1))$ with $(\dd f_1)\sst\C\{x\}$ the Jacobian ideal of $f_1$. (We do not have a problem on the difference between local and global deformations, see for instance \cite[Rem.\,1.3]{sem}.) There is the filtration $V$ on $\C\{x\}$ such that $V^{\al}\C\{x\}$ is generated by monomials whose {\it shifted weighted degree\1} is at least $\al$, where ``shifted" means that the weighted degree is shifted by the minimal spectral number $\alt_f\eq\alt_{f_1}\defs\msum_{i=1}^n\,w_i$ with $w_i$ the weights of variables $x_i$, see Remarks~\hl{R1.1}{1.1} and \hl{R1.2d}{1.2d} below. Put $\rho_f\defs\min\{\be\sgt1\,|\,f_{\be}\nes0\}$.
\sk
Calculating the Gauss-Manin connection in the {\it semi-weighted-homogeneous\1} case, we get the following as a corollary of Theorem~\hl{T1}{1}.
\par\htt{T2}{}\msn
{\bf Theorem~2.} {\it {\rm(i)} Assume $\al\slt\alt_f\pl\rho_f\mi 1$. Then $\nut_{f,\al}\eq\nut_{f_1,\al}$, and hence
\htt{6}{}
$$\ell_{\D_{X,0}}(\D_{X,0}f^{-\al})=\ell_{\D_{X,0}}(\D_{X,0}\1f_1^{-\al}).
\leqno(6)$$
\par\nin \nin
{\rm(ii)} Assume there is $h\ins V^{\al-r\rho_f}\C\{x\}$ with $\al\ins\Q_{>0}$, $r\ins\Z_{>0}$ such that the image of $f_{\!\rho_f}^{\1 r}h$ in $\Gr_V^{\al}\bl(\C\{x\}/(\dd f)\br)$ does not vanish. Then we have the strict inequality
\htt{7}{}
$$\ell_{\D_{X,0}}(\D_{X,0}f^{-\al})>\ell_{\D_{X,0}}(\D_{X,0}\1f_1^{-\al}),
\leqno(7)$$
\par\nin and in the case $h\eq 1$ with $r\eq1$, we have the equality}
\htt{8}{}
$$\ell_{\D_{X,0}}(\D_{X,0}f^{-\al})=\ell_{\D_{X,0}}(\D_{X,0}\1f_1^{-\al})\pl 1,
\leqno(8)$$
\par\nin \ms
This gives certain sufficient conditions for the conjecture \cite[Conj.\,1.7]{BiSc} to hold or to fail, see also \cite{MO}. Note that $n_{f,\al}\eq n_{f_1,\al}$ (see Remark~\hl{R1.2b}{1.2b} below), $\rf\eq r\!_{f_1}$, and $\nut_{f_1,\al}\eq\msum_{j\in\N}\,n_{f,\al-j}$, see (\hl{4}{4}). In the case $f_1\eq x^4{+}y^4{+}z^{4}$ for instance, setting $f_{>1}\eq x^2y^2z^2$ and $x^2y^2z$, the assumptions of (i) and (ii) for $r\eq1$ are satisfied with $\al\eq 1$, $\alt_f\eq\tfrac{3}{4}$, and $\rho_f\eq\tfrac{3}{2}$ and $\tfrac{5}{4}$ respectively (compare with \cite{MO}). More generally, we have the following.
\par\htt{C1}{}\msn
{\bf Corollary~1.} {\it For any homogeneous polynomial $f_1$ in $n$ variables of degree $d\sgt n\gess 3$ having an isolated singularity at $0$, there is always a monomial $f_{>1}$ such that the hypotheses of Theorem~{\rm\hl{T2}{2}}\,{\rm(i)} and {\rm(ii)} with $r\eq 1$ are satisfied for $f\defs f_1\pl f_{>1}$, where $\al\eq 1$, $\alt_f\eq\tfrac{n}{d}$, and the usual degree of $f_{>1}$ is respectively $\,2d{-}n{+}1$ and $\,2d{-}n$.}
\ms
Indeed, for any non-negative integer $k\less n(d{-}2)$, there is a monomial $g_k$ of degree $k$ which does not vanish in the Jacobian ring $\C\{x\}/(\dd f_1)$ using the symmetry of spectral numbers, see (\hl{1.2.3}{1.2.3}) below. We apply this to the case $k\eq 2d{-}n{+}1$ and $2d{-}n$. (Note that $2d{-}n\slt n(d{-}2)$, since $(n{-}2)(d{-}1)\sgt2$.) Corollary~\hl{C1}{1} then follows.
\ms
It does not seem easy to determine the length $\ell_{\D_{X,0}}(\D_{X,0}f^{-\al})$ in the case where the assumption of Theorem~\hl{T2}{2}\,(ii) is satisfied (except the case $r\eq 1$ and $\al\eq\alt_f\pl\rho_f\mi 1$) even though the nilpotent operator $N$ on the right-hand side of (\hl{3}{3}) vanishes in the semi-weighted-homogeneous case. This does not seem trivial even in the case $n_{f,\al}\eq1$ or 0 ($\forall\,\al)$, for instance, if $f_1\eq\msum_{i=1}^n\,x_i^{e_i}$ with $e_i$ mutually prime, where the problem is equivalent to the determination of $b_f(s)$, see \cite{Kat1}, \cite{Kat2}, \cite{Ca}, \cite{sem}, and Section~\hl{S3}{3} below. Note also that the hypothesis of Theorem~\hl{T2}{2}\,(ii) is not a {\it necessary condition\1} for the inequality~(\hl{7}{7}) to hold.
\sk
In Section 1 we review some basics of Brieskorn lattices, spectral numbers, BS polynomials for isolated singularities together with a formula for $\D_Xf^{-\al}/\D_Xf^{-\al+1}$ in \cite[Thm.\,1]{rp}. In Section~2 we prove the main theorems applying the assertions in the previous section. In Section 3 we give some examples and remarks related to the main theorems. In Section~4 we explain a simple proof of a generalization of the last assertion in \cite[Cor.\,3.0.1]{Bit} to the non-integral exponent case implying a simpler proof of Theorem~\hl{T1}{1}.
\sk
This work was partially supported by JSPS Kakenhi 15K04816.
\bs\bs
\vbox{\centerline{\bf 1. Preliminaries}
\bsn
In this section we review some basics of Brieskorn lattices, spectral numbers, BS polynomials for isolated singularities, and a formula for $\D_Xf^{-\al}/\D_Xf^{-\al+1}$ in \cite[Thm.\,1]{rp}.}
\par\htt{1.1}{}\msn
{\bf 1.1.~Brieskorn lattices.} Let $f\ins\C\{x\}$ be a convergent power series of $n$ variables having an isolated singularity at 0, where $f(0)\eq0$ and $n\gess 2$. Set $(X,0):=(\C^n,0)$. The {\it Brieskorn lattice\1} \cite{Br} is defined by
\htt{1.1.1}{}
$$\Hppf:=\Om_{X,0}^n/\ddd f{\wedge}\1\ddd\Om_{X,0}^{n-2}.
\leqno(1.1.1)$$
\par\nin This is a free module of rank $\mu_f$ (with $\mu_f$ the Milnor number) over $\C\{t\}$ and $\C\{\!\{\dd_t^{-1}\}\!\}$, see for instance \cite[1.8]{bl}. The actions of $t$ and $\dd_t^{-1}$ are defined respectively by multiplication by $f$ and
\htt{1.1.2}{}
$$\dd_t^{-1}[\om]\eq[\ddd f{\wedge}\eta]\,\,\,\h{if}\,\,\,\ddd\eta\eq\om\,\,\,(\om\ins\Om_{X,0}^n,\,\eta\ins\Om_{X,0}^{n-1}).
\leqno(1.1.2)$$
\par\nin Here the well-definedness follows from the Poincar\'e lemma. Set
\htt{1.1.3}{}
$$\G_f:=\Hppf[\dd_t],
\leqno(1.1.3)$$
\par\nin which is the localization of $\Hppf$ by the action of $\dd_t^{-1}$. This is a regular holonomic $\D_{X,0}$-module with quasi-unipotent monodromy, and is called the {\it Gauss-Manin system.} It has the $V$-filtration of Kashiwara and Malgrange {\it indexed by\1} $\Q$ so that $\dd_tt{-}\al$ is nilpotent on each graded piece $\G_f^{(\al)}\defs\Gr_V^{\al}\G_f$, see for instance \cite{bl}, \cite{mic}.
\sk
Since $\G_f$ is a regular holonomic $\D_{\C,0}$-module, there are canonical isomorphisms
\htt{1.1.4}{}
$$\G_f=\widehat{\mopl}_{\al\in\Q}\,\G_f^{(\al)},\q\Htppf=\widehat{\mopl}_{\al\in\Q}\,\Htppf{}^{(\al)}.
\leqno(1.1.4)$$
\par\nin Here $\widehat{\mopl}$ is the completion of the direct sum by an appropriate topology (similar to the case of $\C\{t\}\,{\supset}\,\C[t])$, and $\G_f^{(\al)}$ is identified with the kernel of $(\dd_tt\mi\al)^k$ in $\G_f$ with $k\eq\mu_f$ (similarly for $\Htppf{}^{(\al)}$). This is well known for $\G_f$ using a classification of indecomposable regular holonomic $\D_{\C,0}$-modules or the combinatorial description of regular holonomic $\D_{\C,0}$-modules, see for instance \cite[1.3]{bl}, \cite[Rem.\,1.3b]{rh}. For $\Htppf$, there is a decomposition of $\Htppf/V^{\be}\!\Htppf$ by the eigenvalues of $\dd_tt$ for $\be\gg0$, see for instance \cite[Rem.\,A.7c]{rh}. This decomposition is compatible with the projections $\Htppf/V^{\be'}\!\Htppf\onto\Htppf/V^{\be}\!\Htppf$ for $\be\slt\be'$ by the uniqueness of the decomposition. We then get the decomposition of the $V$-adic completion $(\Htppf)^{\wedge}$ of $\Htppf$. This is compatible with the first isomorphism of (\hl{1.1.4}{1.1.4}) applying the same argument with $\Htppf$ replaced by $V^{\al}\G_f$ and taking the inductive limit for $\al$. We have
$$\Htppf/V^{\al'}\Htppf=\mprod_{\al<\al'}\,\Htppf{}^{(\al)},\q\h{hence}\q(\Htppf)^{\wedge}=\mprod_{\al\in\Q}\,\Htppf{}^{(\al)}.$$
\par\nin Indeed, the projection $\rho_{\al}\cols V^{\al}\Htppf\onto\Gr_V^{\al}\Htppf$ splits modulo $V^{\be}\G_f$ in a compatible way with the one for $\G_f$ (hence its $\gamma$-components vanish for $\gamma\ins(\al,\be)$, see (\hl{1.1.5}{1.1.5}) below) and $V^{\be}\G_f\sst\Htppf$ for $\be\,{\gg}\,0$. We then get the canonical splitting of $\rho_{\al}$ whose image is contained in $\G_f^{(\al)}\sst\G_f$ (which is compatible with the one for $\G_f$). This implies the assertion using the decomposition in $\G_f$. Here we can show that $(\Htppf)^{\wedge}\cap\G_f\eq\Htppf$ employing a commutative diagram of short exact sequences and the Mittag-Leffler condition \cite{Gro}, \cite{Ha} (here $\G_f/\Htppf$ does not change by the $V$-adic completion), see also \cite{exp}.
\sk
For $\xi\ins\G_f$, we thus get the {\it asymptotic expansion\1}
\htt{1.1.5}{}
$$\xi=\msum_{\al\in\Q}\,\xi_{\al}\q\h{with}\q\xi_{\al}\ins\G_f^{(\al)}\,\,(\al\ins\Q),
\leqno(1.1.5)$$
\par\nin where $\xi_{\al}$ is called the $\al$-component of $\xi$. Set
\htt{1.1.6}{}
$$\alt(\xi):=\min\{\al\in\Q\mid\xi_{\al}\nes 0\},\q\Gr^V\!\1\xi:=\xi_{\alt(\xi)}.
\leqno(1.1.6)$$
\par\nin \par\htt{R1.1}{}\msn
{\bf Remark 1.1.} In the {\it weighted homogeneous\1} case, we have the Euler field
$$\xi\defs\msum_{i=1}^n\,w_ix_i\dd_{x_i}\q\h{with}\q\xi(f)\eq f.$$
\par\nin Using the {\it interior product\1} $\iota_{\xi}$ and the {\it Lie derivation\1} $L_{\xi}$, it is easy to see that $\Hppf$ is stable by the action of $\dd_tt$, and the $V$-filtration on it is induced by the {\it weighted degree\1} of differential forms, which can be given by the Lie derivation $L_{\xi}$, where $\deg_wx_i\eq \deg_w\ddd x_i\eq w_i$, see also \cite{exp}, \cite[\S 1.4]{BaSa}. Indeed, this follows from the well-known relations
\htt{1.1.7}{}
$$\ddd\1\iota_{\xi}\om=L_{\xi}\1\om,\q\ddd f{\wedge}\iota_{\xi}\1\om=f\om\q(\om\ins\Om_X^n).
\leqno(1.1.7)$$
\par\nin \par\htt{1.2}{}\msn
{\bf 1.2.~Spectral numbers.} The {\it spectrum\1} ${\rm Sp}_f(t)=\msum_{\al\in\Q}\,n_{f,\al}t^{\al}$ was first introduced in \cite{St} for isolated singularities. (There was some confusion about the monodromy, see \cite{DS2}.) The multiplicities $n_{f,\al}$ can be defined by
\htt{1.2.1}{}
$$\aligned n_{f,\al}:={}&\dim_{\C}\Gr^p_FH^{n-1}(F_{\!f},\C)_{\la}\\ ={}&\dim_{\C}\Gr_V^{\al}(\Hppf\!/\dd_t^{-1}\!\Hppf)\\ ={}&\dim_{\C}\Gr_V^{\al}\bl(\C\{x\}/(\dd f)\br),\endaligned
\leqno(1.2.1)$$
\par\nin with $p\eq[n{-}\al]$, $\la\eq e^{-2\pi i\al}$, where $H^{n-1}(F_{\!f},\C)_{\la}$ denotes the generalized $\la$-eigenspace of the monodromy on the vanishing cohomology with $F_{\!f}$ the Milnor fiber. The first equality follows for instance from \cite{SS}, \cite{Va1} (see also \cite{bl},\,\cite[\S 1]{JKSY2}), and the second from the isomorphism
\htt{1.2.2}{}
$$\Hppf\!/\dd_t^{-1}\!\Hppf\cong\C\{x\}/(\dd f),
\leqno(1.2.2)$$
\par\nin trivializing $\Om_{X,0}^n$, where $(\dd f)$ denotes the Jacobian ideal of $f$. We call $\al$ a {\it spectral number\1} of $f$ if $n_{f,\al}\nes0$. It is well known (see for instance \cite{St}, \cite{bl}) that there is the symmetry
\htt{1.2.3}{}
$$n_{f,\al}=n_{f,n-\al}\q(\al\ins\Q).
\leqno(1.2.3)$$
\par\nin \par\htt{R1.2a}{}\msn
{\bf Remark 1.2a.} The equality~(\hl{2}{2}) in the introduction follows from (\hl{1.2.1}{1.2.1}), since we can replace $\dd_t^{-1}\!\Hppf$ in the second term by $t\1\Hppf$ passing to the graded pieces of the monodromy filtration $W$ associated with the nilpotent operator $N\defs\dd_tt{-}\al$ on $\Gr_V^{\al}\G_f$ as in \cite{Va1}. Note that
\htt{1.2.4}{}
$$\dd_t^{-1}\!\Hppf\ne t\1\Hppf\q\h{(although}\,\,\,\,\dd_t^{-1}\!\Htppf\eq t\1\Htppf),
\leqno(1.2.4)$$
\par\nin unless $f$ is weighted homogeneous. (This is related to the last part of Remark~\hl{R1.3a}{1.3a} below.)
\par\htt{R1.2b}{}\msn
{\bf Remark 1.2b.} It is well known (see for instance \cite{Va2}) that the spectrum ${\rm Sp}_f(t)$ is invariant under a $\mu$-constant deformation of $f$ (e.g., under a semi-weighted-homogeneous deformation of a weighted homogeneous polynomial with an isolated singularity).
\par\htt{R1.2c}{}\msn
{\bf Remark 1.2c.} In the weighted homogeneous case, the $V$-filtration on the Jacobian ring $\C\{x\}/(\dd f)$ can be given by using the {\it shifted weighted degree\1} as in Remark~\hl{R1.1}{1.1}. Here we have the shift coming from the trivialization of the highest differential forms $\Om_{X,0}^n$, which is given by the sum of the weights of variables $\msum_{i=1}^n\,w_i$ in the notation of Theorem~\hl{T2}{2}.
\par\htt{R1.2d}{}\msn
{\bf Remark 1.2d.} Let $f\eq f_1\pl f_{>1}$ be a semi-weighted-homogeneous deformation of a weighted homogeneous polynomial $f_1$ having an isolated singularity as in Theorem~\hl{T2}{2}. Let $V$ be the decreasing filtration on $\C\{x\}$ defined by the condition that the {\it shifted\1} weighted degree is at least $\al$ as in Remark~\hl{R1.2c}{1.2c} just above. This induces the $V$-filtration on the Gauss-Manin system $\G_f$ (see Remark~\hl{R1.1}{1.1} for $\G_{f_1}$). This can be verified for instance by using the $V$-adic completion together with the Mittag-Leffler condition (see \cite[Prop.\,13.2.3]{Gro}, \cite{Ha}) as in \cite{exp}. Here we have the canonical isomorphism
\htt{1.2.5}{}
$$\Gr_V^{\ssb}\bl(\C\{x\}/(\dd f)\br)=\Gr_V^{\ssb}\bl(\C\{x\}/(\dd f_1)\br),
\leqno(1.2.5)$$
\par\nin where the weighted degree is shifted by $\msum_{i=1}^n\,w_i$ as in Remark~\hl{R1.2c}{1.2c}.
This isomorphism can be verified by using the filtered shifted Koszul complexes
$$\bl((\Om_{X,0}^{\ssb},V),\ddd f\wedge\br)[n],\q\bl((\Om_{X,0}^{\ssb},V),\ddd f_1\wedge\br)[n],$$
\par\nin which give {\it filtered\1} free resolutions of the Jacobian rings $\C\{x\}/(\dd f)$, $\C\{x\}/(\dd f_1)$ respectively.
\par\htt{1.3}{}\msn
{\bf 1.3.~BS polynomials.} Let $f$ and $\Hppf,\G_f$ be as in \hl{1.1}{1.1}. We have the {\it saturation\1} $\Htppf\sst\G_f$ defined as in (\hl{1}{1}) in the introduction. Let $b_f(s)$ be the BS polynomial of $f$. Set
\htt{1.3.1}{}
$$\bt_f(s):=b_f(s)/(s{+}1).
\leqno(1.3.1)$$
\par\nin This is called the {\it reduced\1} (or {\it microlocal\1}) BS polynomial, see \cite{mic}. Note that the roots of $\bt_f(s)$ are strictly negative rational numbers, see \cite{Ka}. By B.~Malgrange \cite{Ma}, we have the following equality as is explained in the introduction:
\htt{1.3.2}{}
$$\bt_f(s)={\rm Mini.\,Poly}\bl({-}\dd_tt\ins{\rm End}_{\C}(\Htppf/t\Htppf)\br).
\leqno(1.3.2)$$
\par\nin Comparing this with (\hl{2}{2}) in the introduction, we see that the roots of $\bt_f(s)$ coincides with the spectral numbers up to sign in the weighted homogenous case (where $\Htppf\eq\Hppf$).
\par\htt{R1.3a}{}\msn
{\bf Remark 1.3a.} Let $\al_{f,k}$ ($k\ins[1,\mu_f]$) be the spectral numbers of $f$ written with multiplicities. It is quite well known that there are {\it non-negative\1} integers $r_k$ ($k\ins[1,\mu_f]$) such that the $\al_{f,k}\mi r_k$ ($k\ins[1,\mu_f]$) are the roots of $\bt_f(s)$ up to sign forgetting the multiplicities, and
\htt{1.3.3}{}
$$\min\{\al_{f,k}\}_{k\in[1,\mu_f]}=\min\{\al_{f,k}\mi r_k\}_{k\in[1,\mu_f]}.
\leqno(1.3.3)$$
\par\nin These easily follow from \cite{Ma}, \cite{SS}, \cite{Va1}. Here we define the filtration $\Ft^{\ssb}\,{\supset}\,F^{\ssb}$ on $H^{n-1}(F_{\!f},\C)_{\la}$ replacing $\Hppf$ with $\Htppf$ in \cite{SS}, see also \cite[(2.6.3)]{bl}. Then the $r_k$ can be obtained by considering the meaning of the bigraded pieces $\Gr_{\Ft}^{\ssb}\Gr_F^{\ssb}H^{n-1}(F_{\!f},\C)_{\la}$.
\sk
It is also well known that a convergent power series $f$ with an isolated singularity is weighted homogeneous after some coordinate change if and only if the above $r_k$ vanish for any $k$ (using \cite{SaK} together with a formalism of primitive forms if necessary, see for instance a remark after Proposition in the introduction of \cite{bl}). Recall that a convergent power series with an isolated singularity has finite determinacy, see for instance \cite{GLS}.
\par\htt{R1.3b}{}\msn
{\bf Remark 1.3b.} It is well known that any (small) $\mu$-constant deformation of a weighted homogeneous polynomial having an isolated singularity can be given by a semi-weighted-homogeneous deformation, see \cite{Va3}. (This implies the coincidence of {\it modality\1} and {\it inner modality\1} for weighted homogeneous polynomials with isolated singularities.) There is a stratification of the parameter space of the miniversal $\mu$-constant deformation such that $b_f(s)$ is constant on each stratum. However, it is not necessarily easy to determine this explicitly. Using for instance ``bernstein" in \cite[gmssing.lib]{Sing}, one may verify the computations of examples in \cite{Kat1}, \cite{Kat2}, \cite{Ca} at {\it typical points\1} of strata. The former looks compatible with this. The calculation for the latter seems rather nontrivial, see also \cite{sem} and Remark~\hl{R3.5b}{3.5b} below.
\par\htt{1.4}{}\msn
{\bf 1.4.~Description of the quotient $\D_Xf^{-\al}/\D_Xf^{-\al+1}\,$.} For $X,f$ as in \hl{1.1}{1.1}, we denote by $i_f\cols X\,{\into}\,Y\defs X{\times}\C$ the graph embedding for $f$. Set
$$\B_f:=\bl((i_f)^{\D}_*\OO_X\br){}_0\,\,\bl(=\OO_{X,0}[\dd_t]\delta(t{-}f)\br),$$
\par\nin where $(i_f)^{\D}_*$ is the direct image as $\D$-module, see for instance\cite[\S 1]{rp}. It is a regular holonomic $\D_{Y,0}$-module, and has the $V$-filtration of Kashiwara and Malgrange indexed by $\Q$ so that $N\defs\dd_tt{-}\al$ is nilpotent on $\Gr_V^{\al}$, where $t$ is the coordinate of $\C$. Put
$$\Mf^{(\al)}:=\Gr_V^{\al}\B_f\q(\al\ins\Q).$$
\par\nin These are regular holonomic $\D_{X,0}$-modules (corresponding to the $\la$-eigenspace of the nearby cycle complex $\psi_{f,\la}\C_X[n{-}1]$ with $\la\defs e^{-2\pi i\al}$). Using the above property of $\dd_tt$ on $\Gr_V^{\al}$, we can deduce the isomorphisms
\htt{1.4.1}{}
$$t\cols\Mf^{(\al)}\!\simto\!\Mf^{(\al+1)},\,\,\,\,\dd_t\cols\Mf^{(\al+1)}\!\simto\!\Mf^{(\al)}\,\,\,\,\,\,(\al\nes0).
\leqno(1.4.1)$$
\par\nin \sk
As is well known (see for instance \cite{Ma}), there is a natural inclusion
$$\Nf:=\D_{X,0}[s]f^s\into\B_f,$$
\par\nin where $s$ and $f^s$ are identified with $-\dd_tt$ and $\delta(t{-}f)$ respectively, and the BS polynomial $b_f(s)$ coincides with the minimal polynomial of the action of $s\eq{-}\dd_tt$ on the quotient regular holonomic $\D_{X,0}$-module
\htt{1.4.2}{}
$$\Nf/t\Nf,
\leqno(1.4.2)$$
\par\nin where the action of $t$ on $\Nf$ is defined by $s\mapsto s{+}1$, see \cite{Ka}. Note that the quotient filtration $V$ on $\Nf/t\Nf$ is a finite filtration (since $\Nf/t\Nf$ is holonomic), and $-\al$ is a root of $b_f(s)$ if and only if $\Gr_V^{\al}(\Nf/t\Nf)\ne0$. Set
$$\Nf{}^{(\al)}:=\Gr_V^{\al}\!\Nf\,\subset\,\Mf^{(\al)}\q(\al\ins\Q).$$
\par\nin \sk
We define the filtration $G$ on $\Mf^{(\be)}$ for $\be\ins(0,1]$ so that
\htt{1.4.3}{}
$$t^j:G_j\Mf^{(\be)}\simto\Nf{}^{(\be+j)}\,\,\bl(\subset\Mf^{(\be+j)}\br)\q(j\ins\N),
\leqno(1.4.3)$$
\par\nin where $G_j\Mf^{(\be)}\eq0$ for $j\slt0$. This definition coincides with the one in \cite{rp}, since there is the inclusion
\htt{1.4.4}{}
$$\Nf\sst V^{>0}\B_f,
\leqno(1.4.4)$$
\par\nin by negativity of the roots of $b_f(s)$, see \cite{Ka}. We have also the isomorphisms
\htt{1.4.5}{}
$$\dd_t^j:\Nf{}^{(\be+j)}\simto G_j\Mf^{(\be)}\q(j\ins\N),
\leqno(1.4.5)$$
\par\nin since the holonomic $\D_{X,0}$-submodule $\Nf{}^{(\be+j)}$ is {\it stable\1} by the {\it automorphism\1} $t^j\dd_t^j$ of $\Mf^{(\be+j)}$ (and an injective endomorphism of a holonomic $\D$-module is an automorphism). Indeed, we have
\htt{1.4.6}{}
$$t^j\dd_t^j\eq\mprod_{k=1}^j\,(\dd_tt{-}k)\eq\mprod_{k=1}^j\,(N{+}\be{+}j{-}k)\q\h{on}\q\Mf^{(\be+j)},
\leqno(1.4.6)$$
\par\nin and $\Nf{}^{(\be+j)}\sst\Mf^{(\be+j)}$ is stable by the action of $N\eq{-}s{-}\be{-}j$ on $\Mf^{(\be+j)}$.
\sk
We have the following
\par\htt{T1.4}{}\msn
{\bf Theorem~1.4} (\cite[Thm\,1]{rp}). {\it There are isomorphisms of $\D_{X,0}$-modules}
\htt{1.4.7}{}
$$\D_{X,0}f^{-\be-j}\!/\D_{X,0}f^{-\be-j+1}=\Gr^G_j(\Mf^{(\be)}\!/N\Mf^{(\be)})\q(\be\ins(0,1],j\ins\N).
\leqno(1.4.7)$$
\par\nin \ms
Here $G$ denotes also the quotient filtration on $\Mf^{(\be)}\!/N\Mf^{(\be)}$ so that we have the {\it strict surjection\1}
$$\bl(\Mf^{(\be)},G\br)\onto\bl(\Mf^{(\be)}\!/N\Mf^{(\be)},G\br).$$
\par\nin Note that $N$ is not necessarily {\it strictly compatible\1} with $G$ in general, see \cite{rp}.
\sk
We need also the following for the proof of Theorem~\hl{T1}{1}.
\par\htt{P1.4}{}\msn
{\bf Proposition~1.4} (\cite[Prop.\,1]{rp}). {\it Setting
$$\Mtf^{\,(1)}:=\Mf^{(1)}/{\rm Ker}\,N\,\,\bl(=\Mf^{(0)}\br),$$
\par\nin we have the isomorphisms
\htt{1.4.8}{}
$$\Gr^G_j\Mf^{(1)}\simto\Gr^G_j\Mtf^{\,(1)},\q\q\h{}
\leqno(1.4.8)$$
\par\nin \vskip-6mm\htt{1.4.9}{}
$$\h{}\q\q\Gr^G_j\bl(\Mf^{(1)}/N\Mf^{(1)}\br)\simto\Gr^G_j\bl(\Mtf^{\,(1)}/N\Mtf^{\,(1)}\br)\q(j\gess 1),
\leqno(1.4.9)$$
\par\nin and the short exact sequences of regular holonomic $\D_{X,0}$-modules}
\htt{1.4.10}{}
$$0\to \M_Z\to\Gr^G_0\Mf^{(1)}\to\Gr^G_0\Mtf^{\,(1)}\to 0,
\leqno(1.4.10)$$
\par\nin \vskip-6mm\htt{1.4.11}{}
$$0\to \M_Z^{\rm IC}\to\Gr^G_0(\Mf^{(1)}/N\Mf^{(1)})\to\Gr^G_0\bl(\Mtf^{\,(1)}/N\Mtf^{\,(1)}\br)\to 0.
\leqno(1.4.11)$$
\par\nin \sk
Here $\M_Z$, $\M_Z^{\rm IC}$ are regular holonomic $\D_{X,0}$-modules such that $\DR_X(\M_Z)\eq\C_Z[n{-}1]$ and $\DR_X(\M_Z^{\rm IC})\eq\h{\rm IC}_Z\C$ with $\h{\rm IC}_Z\C$ the intersection complex of $Z:=f^{-1}(0)\sst X$. These $\D_{X,0}$-modules can be extended to $\D_X$-modules by shrinking $X$ sufficiently.
\par\htt{C1.4}{}\msn
{\bf Corollary~1.4.} {\it For $\be\ins(0,1]$, we have}
\htt{1.4.12}{}
$${\rm Supp}\,\Gr^G_j\Mf^{(\be)}\subset\{0\}\q\h{if}\q(\be,j)\ne(1,0).
\leqno(1.4.12)$$
\par\nin \bs\bs
\vbox{\centerline{\bf 2. Proof of the main theorems}
\bsn
In this section we prove the main theorems applying the assertions in the previous section.}
\par\htt{2.1}{}\msn
{\bf 2.1.~Proof of Theorem~\hl{T1}{1}.} It is well known that we have a canonical isomorphism of regular holonomic $\D_{\C,0}$-modules
\htt{2.1.1}{}
$$\Hc^0\DR_X(\B_f)=\G_f,
\leqno(2.1.1)$$
\par\nin see for instance \cite{fil}, \cite{SS}, \cite{exp}. Moreover, the $V$-filtration on the complex $\DR_X(\B_f)$ is {\it strict,} and induces the $V$-filtration on $\G_f$. This can be verified by using for instance an argument similar to the proof of \cite[Prop.\,3.4.8]{mhp}, see also \cite[\S 4.11]{DS1}, \cite[\S 1.4]{JKSY2}.
\sk
We denote $\Hc^j\DR_X$ by $\Hc^j_{\DR}$ to simplify the notation. In the notation of \hl{1.4}{1.4}, set
$$\Cc_f:=\B_f/\Nf.$$
\par\nin From (\hl{1.4.10}{1.4.10}) we can deduce that
\htt{2.1.2}{}
$${\rm Supp}\,\Cc_f\subset\{0\},
\leqno(2.1.2)$$
\par\nin since ${\rm Supp}\,\Mf^{(\al)}\sst\{0\}$ ($\forall\,\al\slt1$). Hence $\Cc_f$ is isomorphic to an infinite direct sum of (copies of) $\B_{\{0\}}:=\C[\dd_{x_1},\dots,\dd_{x_n}]$ by a special case of Kashiwara's equivalence as in \cite{Ma}, see also \cite[A.7]{rh}.
We then get for $j\ins\Z$ the commutative diagram of short exact sequences
\htt{2.1.3}{}
$$\begin{array}{ccccccccccc}&&0&&0&&0\\
&&\downarrow&&\downarrow&&\downarrow&\raise-2mm\h{}\\
0&\to&\Hc^j_{\DR}V^{\al}\!\Nf&\to&\Hc^j_{\DR}V^{\al}\B_f&\to&\Hc^j_{\DR}V^{\al}\Cc_f&\to&0\\
&&\downarrow&&\downarrow&&\downarrow&\raise4mm\h{}\raise-2mm\h{}\\
0&\to&\Hc^j_{\DR}\Nf&\to&\Hc^j_{\DR}\B_f&\to&\Hc^j_{\DR}\Cc_f&\to&0\\
&&\downarrow&&\downarrow&&\downarrow&\raise4mm\h{}\raise-2mm\h{}\\
0&\to&\Hc^j_{\DR}\Nf/V^{\al}\!\Nf&\to&\Hc^j_{\DR}\B_f/V^{\al}\B_f&\to&\Hc^j_{\DR}\Cc_f/V^{\al}\Cc_f&\to&0\\
&&\downarrow&&\downarrow&&\downarrow&\raise4mm\h{}\raise-1mm\h{}\\&&0&&0&&0\\
\end{array}
\leqno(2.1.3)$$
\par\nin since $\Hc^j_{\DR}\Cc_f\eq0$ ($j\nes0$) by (\hl{2.1.2}{2.1.2}) (similarly with $\Cc_f$ replaced by $V^{\al}\Cc_f$, $\Cc_f/V^{\al}\Cc_f$) and $\DR_X(\Nf)^j\eq0$ ($j\sgt0$). Indeed, these imply the exactness of the rows and the right column, and that of the middle column follows from the strictness of the filtration $V$ on the complex $\DR_X(\B_f)$ explained above. We then get the exactness of the left column by the nine lemma.
\sk
The last exactness means the strictness of the filtration $V$ on the complex $\DR_X(\Nf)$, and the diagram implies the {\it strict injectivity\1} of the canonical morphism
\htt{2.1.4}{}
$$(\Hc^0_{\DR}\Nf,V)\into(\Hc^0_{\DR}\B_f,V)\eq(\G_f,V),
\leqno(2.1.4)$$
\par\nin where the last isomorphism follows from (\hl{2.1.1}{2.1.1}) and a remark after it.
Then $\Hc^0_{\DR}\Nf$ is identified with the image of $\Om_{X,0}^n[s]$ in $\G_f$, since the de Rham complex is the Koszul complex for the actions of the $\dd_{x_i}$ fixing the coordinates $x_1,\dots,x_n$. Moreover the image of $\Om_{X,0}^n[s]$ is $\Htppf$, since that of $\Om_{X,0}^n$ is $\Hppf$. We thus get the canonical isomorphisms
\htt{2.1.5}{}
$$(\Htppf,V)=(\Hc^0_{\DR}\Nf,V).
\leqno(2.1.5)$$
\par\nin \sk
We have also the following isomorphism for $\al\ins\Q$ using a similar diagram whose columns are the short exact sequence $0\tos V^{>\al}\tos V^{\al}\tos\Gr_V^{\al}\tos 0$\,:
\htt{2.1.6}{}
$$\begin{array}{ccccccc}\Htppf{}^{(\al)}&=&\Hc^0_{\DR}\Nf{}^{(\al)}\\ \rotatebox{90}{$\supset$}&\raise5mm\h{}&\rotatebox{90}{$\supset$}\\ \,\,\G_f^{(\al)}&=&\,\,\Hc^0_{\DR}\Mf{}^{(\al)}.\end{array}
\leqno(2.1.6)$$
\par\nin Here the last isomorphism follows from a remark after (\hl{2.1.1}{2.1.1}). (It is not easy to prove (\hl{2.1.6}{2.1.6}) without using the above argument when $n\eq2$.)
\sk
In the case $\al\,{\notin}\,\Z$, the equality~(\hl{5}{5}) then follows from Theorem~\hl{T1.4}{1.4} and Corollary~\hl{C1.4}{1.4}. Indeed, the functor $\Hc^0_{\DR}$ induces an equivalence between the {\it abelian category\1} of regular holonomic $\D_X$-modules supported on $\{0\}$ and that of finite dimensional $\C$-vector spaces (in particular, it commutes with the image of morphisms). Note that the {\it Verdier-type extension theorem\1} is used in the proof of Theorem~\hl{T1.4}{1.4}, where the image of $N$ corresponds to the {\it minimal\1} extension, which has length 1. This is the reason for which we have to consider module the image of $N$ on the right-hand side of (\hl{1.4.7}{1.4.7}).
\sk
In the case $\al\ins\Z$, we can replace $\Mf^{(1)}$ with $\Mtf^{\,(1)}$ (which corresponds to the vanishing cycle complex $\varphi_{f,1}\C_X[n{-}1]$) using Proposition~\hl{P1.4}{1.4}. Indeed, the image by the functor $\Hc^0_{\DR}$ does not change by this (using (\hl{1.4.8}{1.4.8}), (\hl{1.4.10}{1.4.10})). On the other hand, the length of the $\D_{X,0}$-module on the right-hand side of (\hl{1.4.7}{1.4.7}) decreases by $\rf\delta_{\al,1}$ (using (\hl{1.4.9}{1.4.9}), (\hl{1.4.11}{1.4.11})). The assertion then follows by an argument similar to the case $\al\,{\notin}\,\Z$. This completes the proof of Theorem~\hl{T1}{1}.
\par\htt{2.2}{}\msn
{\bf 2.2.~Proof of Theorem~\hl{T2}{2}.} To apply Theorem~\hl{T1}{1}, we have to study the saturation $\Htppf$. For this we calculate the action of $\dd_tt$ on $[h\ddd x]$, which is the image of $h\ddd x$ in $\G_f$ for $h\ins\C\{x\}$. Here we may assume that $h$ is a polynomial and the {\it weighted degree\1} of $h$ is {\it pure,} see Remark~\hl{R1.1}{1.1}. (This means that $h$ is a weighted homogeneous polynomial such that the weights of variables are constant multiples of those for $f_1$ if the weighted degree of $h$ is normalized to 1.) Let $\alt(h)$ be the {\it shifted weighted degree\1} of $h$, see Remark~\hl{R1.2d}{1.2d}. We see that
\htt{2.2.1}{}
$$\aligned&(\dd_tt{-}\alt(h))[h\ddd x]=\msum_{\be>1}\,c_{h,\be}\1\dd_t[f_{\be}h\ddd x],\\
&\h{that is,}\q(t{-}\alt(h)\dd_t^{-1})[h\ddd x]=\msum_{\be>1}\,c_{h,\be}\1[f_{\be}h\ddd x],\endaligned
\leqno(2.2.1)$$
\par\nin with $c_{h,\be}\ins\C^*$ by an argument similar to \cite{bl}, \cite{rp} using (\hl{1.1.2}{1.1.2}), (\hl{1.1.7}{1.1.7}) with $\xi$ the Euler field satisfying $\xi(f_1)\eq f_1$. (Actually one can verify that $c_{h,\be}\eq1{-}\be$.) This can be used to get the {\it asymptotic expansion\1} of $[h\ddd x]\ins\G_f$ inductively. Note that the action of $\dd_tt$ on $\G_f$ is {\it semisimple\1} (since the monodromy $T$ is semisimple). This follows also from (\hl{2.2.1}{2.2.1}).
\sk
We get from (\hl{2.2.1}{2.2.1}) that
\htt{2.2.2}{}
$$(\dd_tt{-}\alt(h))[h\ddd x]=0\q{\rm mod}\,\,\,\,\,V^{\alt_f+\rho_f-1}\G_f.
\leqno(2.2.2)$$
\par\nin Combining this with Theorem~\hl{T1}{1}, Remark\,\,\hl{R1.2d}{1.2d} and (\hl{1.2.5}{1.2.5}), and using the decomposition (\hl{1.1.4}{1.1.4}), we get the assertion (i).
\sk
For the assertion (ii), we have to calculate the {\it asymptotic expansion\1} of $[h\ddd x]$ for $h$ in (ii), see (\hl{1.1.5}{1.1.5}). Set
$$f_{\!\B}:=\mprod_{j=1}^k\,f_{\be_j}\q\h{for}\q\B\eq(\be_1,\dots,\be_k)\ins\Q_{>1}^k.$$
\par\nin Applying (\hl{2.2.1}{2.2.1}) repeatedly, we see that the image of $[h\ddd x]$ in $\G_f/V^{>\al}\G_f$ with $\al$ as in (ii) is contained in the $\C$-vector subspace spanned by
\htt{2.2.3}{}
$$\Gr^V\dd_t^{|\B|}\bl[f_{\!\B}h\ddd x\br]\q(\B\ins\Q_{>1}^k,\,k\ins[0,r]),
\leqno(2.2.3)$$
\par\nin in the notation of (\hl{1.1.6}{1.1.6}). (Note that $N\eq 0$ and the action of $\dd_tt$ is semisimple in the semi-weighted-homogeneous case.) Here $|\B|\eq k$ if $\B\ins\Q_{>1}^k$, and $f_{\!\B}\eq1$ if $|\B|\eq0$. We may assume $k\less r$ and moreover $\B\eq(\rho_f,\dots,\rho_f)$ in the case $k\eq r$, since
\htt{2.2.4}{}
$$\dd_t^{|\B|}\bl[f_{\!\B}h\ddd x\br]\ins V^{>\al}\G_f\,\,\,\h{if}\,\,\,k\sgt r\,\,\,\h{or}\,\,\,k\eq r,\,\B\nes(\rho_f,\dots,\rho_f).
\leqno(2.2.4)$$
\par\nin On the other hand, there is the Hodge filtration $F$ on $\G_f$ with $F_p\G_f:=\dd_t^p\Hppf$ ($p\ins\Z$), and
\htt{2.2.5}{}
$$\dd_t^{|\B|}\bl[f_{\!\B}h\ddd x\br]\ins F_{r-1}\G_f\q\h{if}\q|\B|\slt r.
\leqno(2.2.5)$$
\par\nin We have moreover the canonical isomorphisms
\htt{2.2.6}{}
$$\Gr_V^{\al+r}(\C\{x\}/(\dd f)=\Gr_V^{\al+r}\Gr^F_0\G_f=\Gr_V^{\al}\Gr^F_r\G_f,
\leqno(2.2.6)$$
\par\nin where the last isomorphism is induced by $\dd_t^{\1 r}$. The assertion~(ii) is then reduced the canonical injection
\htt{2.2.7}{}
$$\Gr_V^{\al}\Gr^F_r\G_f\into\G_f/(V^{>\al}\G_f\pl F_{p-1}\G_f),
\leqno(2.2.7)$$
\par\nin which can be verified for instance by using the {\it strict injectivity\1} of the inclusion
\htt{2.2.8}{}
$$(\Gr^F_r\G_f,V)\into(\G_f/F_{r-1}\G_f,V).
\leqno(2.2.8)$$
\par\nin This finishes the proof of Theorem~\hl{T2}{2}.
\par\htt{R2.2}{}\msn
{\bf Remark 2.2.} Assume $f_1\eq\msum_{i=1}^n\,x_i^{e_i}$ and $f_{>1}$ is a monomial $x^a$ with $a\eq(a_1,\dots,a_n)\ins\N^n$. In this case the numbers $r_k$ ($k\ins[1,\mu_f]$) in Remark~\hl{R1.3a}{1.3a} can be calculated rather easily. Indeed, it is quite well known that the {\it spectral numbers\1} of $f$ are given with multiplicities by
\htt{2.2.9}{}
$$S_f=\bl\{\msum_{i=1}^n\,\nu_i/e_i\mid\nu_i\in\Z\,{\cap}\,[1,e_i{-}1]\,\,(i\ins[1,n])\br\}.
\leqno(2.2.9)$$
\par\nin see for instance \cite[\S 1.5]{JKSY2}. For each $\nu\eq(\nu_1,\dots,\nu_n)\ins\Xi\defs\mprod_{i=1}^n\,\Z\,{\cap}\,[1,e_i{-}1]$, set
\htt{2.2.10}{}\htt{2.2.11}{}
$$\Psi(\nu)=\{(\nu'\!,j)\ins\Xi{\times}\Z_{>0}\mid\nu'\pl ja\eq\nu\,\,\,{\rm mod}\,\,\,\mprod_{i=1}^n\,e_i\Z\br\},
\leqno(2.2.10)$$
\par\nin $$\aligned&r(\nu)\eq\max\{r(\nu'\!,j)\}_{(\nu',j)\in\Psi(\nu)}\q\h{with}\\&r(\nu'\!,j)\eq j\mi\msum_{i=1}^n\,[(\nu'_i\pl ja_i)/e_i],\endaligned
\leqno(2.2.11)$$
\par\nin where $r(\nu)\eq0$ if $J(\nu)\eq\emptyset$ (and $[*]$ denotes the integer part, that is, the Gauss symbol). Then these $r(\nu)$ ($\nu\ins\Xi$) give the $r_k$ ($k\ins[1,\mu_k]$) in Remark~\hl{R1.3a}{1.3a}. This follows from the argument in the proof of Theorem~\hl{T2}{2}. (When we calculate the acton of $\dd_t$ on $\Gr^V[x^{\nu}\ddd x]$, we can consider the elements mod $V^{>\al}\G_f$ with $\al\eq\alt(x^{\nu})$ so that the division of $\dd_{x_i}f$ can be replaced by that of $\dd_{x_i}f_1$. This means that the leading term $\Gr^V[x^{\nu}\ddd x]$ in $\G_f$ can be identified with that in $\G_{f_1}$, see also Remark~\hl{R1.2d}{1.2d}. In particular, we have $\Gr^V[x^{\nu}\ddd x]\eq 0$ if $\nu_i{+}1\ins e_i\Z$ for some $i$.)
\bs\bs\htt{S3}{}
\vbox{\centerline{\bf 3. Examples and Remarks}
\bsn
In this section we give some examples and remarks related to the main theorems.}
\par\htt{3.1}{}\msn
{\bf 3.1.~General setting.} Let $f_1=\msum_{i=1}^n\,x_i^{e_i}$ with $e_i$ ($i\ins[1,n]$) integers greater than 2. Let $f_{>1}$ be a linear combinations of monomials of unshifted weighted degree bigger than 1, where $f_1$ has unshifted weighted degree 1. Set $f\eq f_1\pl f_{>1}$. The spectrum of $f$ (and $f_1$) is given by (\hl{2.2.9}{2.2.9}).
\sk
If the $e_i$ are {\it mutually prime,} we have $n_{f,\al}\eq1$ or $0$ ($\forall\,\al$), and moreover the spectral numbers are mutually different modulo $\Z$, that is, the multiplicity of each eigenvalue of the Milnor monodromy is 1. Indeed, for $\nu_i\ins\Z$ ($i\ins[1,n]$), we have
\htt{3.1.1}{}
$$\msum_{i=1}^n\,\nu_i/e_i\ins\Z\iff\nu_i/e_i\ins\Z\,\,\,(\forall\,i\ins[1,n]).
\leqno(3.1.1)$$
\par\nin So the determination of $\Htppf$ is equivalent to that of $b_f(s)$. The latter, however, does not seem quite easy even in the case $f_{>1}$ is a monomial, since it is rather complicated to calculate $[x^{\nu}\ddd x]_{\al}$ $(\al\ins\Q)$ in the notation of (\hl{1.1.5}{1.1.5}) especially when $\nu_i\pl 1\ins e_i\1\Z$, see \hl{3.3}{3.3}--\hl{3.5}{5} below.
\par\htt{3.2}{}\msn
{\bf 3.2.~Example~I.} Let
$$f_1=\msum_{i=1}^6\,x_i^{2a_i},\q f_{>1}=\mprod_{i=1}^6\,x_i^{(a_i-1)/2},$$
\par\nin where the $a_i$ are mutually prime odd positive integers ($i\ins[1,6]$). Then
\htt{3.2.1}{}
$$\rho_f=\msum_{i=1}^6\,\tfrac{a_i-1}{2}\tfrac{1}{2a_i}=\tfrac{3}{2}-\tfrac{1}{4}\1\msum_{i=1}^6\,\tfrac{1}{a_i}>1,
\leqno(3.2.1)$$
\par\nin since
$$\msum_{i=1}^6\,\tfrac{1}{a_i}\les\tfrac{1}{3}{+}\tfrac{1}{5}{+}\tfrac{1}{7}{+}\tfrac{1}{11}{+}\tfrac{1}{13}{+}\tfrac{1}{17}\slt1.$$
\par\nin Moreover we have $\alt(f_{>1}^{\1 2})\eq 3$, and 3 is the unique integral spectral number of $f$. Indeed, there is no $\nu\eq(\nu_i)\in\N^6$ such that $\nu_i\ins[0,2a_i{-}2]$ ($\forall\,i$) and $\alt(x^{\nu})\ins\tfrac{1}{2}\1\Z$ except the case $\nu_i\eq a_i{-}1$ ($\forall\,i$), since the $a_i$ are mutually prime. So this gives an example where the assumption of Theorem~\hl{T1}{2} (ii) is satisfied for $r\eq2$, but not for $r\eq 1$, with $\al\eq1$ fixed. (Here we apply (\hl{2.2.1}{2.2.1}) to $[\ddd x]$ and $[h\ddd x]$.)
\par\htt{3.3}{}\msn
{\bf 3.3.~Example~II.} Let
$$f_1\eq x^6\pl y^5,\q f_{>1}\eq x^2y^4.$$
\par\nin Here $f$ {\it cannot\1} be weighted homogeneous by any coordinate change, although $f_{>1}$ belongs to the Jacobian ideal of $f_1$. Indeed, the {\it Tjurina number\1} of $f$ is 19, which does not coincide with the Milnor number 20, according to Singular \cite{Sing}. By (\hl{2.2.9}{2.2.9}) the set of spectral numbers of $f,f_1$ is given by
\htt{3.3.1}{}
$$S_f=S_{f_1}=\bl\{\tfrac{5i+6j}{30}\mid (i,j)\in[1,5]{\times}[1,6]\br\}.
\leqno(3.3.1)$$
\par\nin This coincides with the set of roots of $\bt_{f_1}(s)$ up to sign, which is denoted by $\Rt_{f_1}$, and that for $\bt_f(s)$ is given by
\htt{3.3.2}{}
$$\Rt_f=S_f\,{\cup}\,\bl\{\tfrac{19}{30}\br\}\stm\bl\{\tfrac{49}{30}\br\}\q\h{with}\q\tfrac{49}{30}\eq\max S_f.
\leqno(3.3.2)$$
\par\nin The latter assertion follows from Remark~\hl{R2.2}{2.2}. It is also possible to reduce it to \cite{SaK} using the last remark in Remark~\hl{R1.3a}{1.3a} and Theorem~\hl{T2}{2}\,(i), since $\alt_f\pl\rho_f\eq\tfrac{45}{30}$ is bigger than the second largest spectral number $\tfrac{44}{30}$, where $\alt_f\eq\tfrac{11}{30}$, $\rho_f\eq\tfrac{34}{30}$. We then get the inequality (\hl{7}{7}) for $\al\eq\tfrac{19}{30}$, which is equal to the maximal spectral number minus 1 (and to $\alt_f\pl2(\rho_f{-}1)$), although the hypothesis of (ii) is not satisfied for $r\eq1,2$.
\par\htt{3.4}{}\msn
{\bf 3.4.~Example~II$'$.} Let
$$f_1=x^6\pl y^5,\q f_{>1}=10\1x^4y^3\pl 5\1x^2y^4.$$
\par\nin We can verify that $f$ is {\it weighted homogenous\1} after a local coordinate change (calculating $(y\pl x^2)^5$). Note that the spectrum of a weighted homogenous isolated singularity depends only on the weights $(w_i)$, see for instance \cite{JKSY2}. Generalizing the calculation in \hl{3.3}{3.3}, we see that some {\it cancellation\1} occurs as in \hl{3.5}{3.5} below. This shows that Theorem~\hl{T2}{2}\,(ii) does not hold without the assumption on the image of $f_{\!\rho_f}^{\1 r}h$ with $r\eq 2$, and we can prove the assertion only for the minimal degree part of $f_{>1}$.
\par\htt{3.5}{}\msn
{\bf 3.5.~Example~III.} Let
$$f_1\eq x^7\pl y^5,\q f_{>1}\eq h_1\pl c\1h_2,\q h_1\eq x^3y^3,\q h_2=x^5y^2,$$
\par\nin where $c\ins\C^*$, see \cite{Kat1}, \cite{sem}. Put $r:=4$. Then
\htt{3.5.1}{}
$$h_2\eq h_1^{\,r}/x^7y^{10},\q\al'(h_1){-}1\eq\tfrac{1}{35},\q\al'(h_2){-}1\eq\tfrac{r}{35},\q\alt_f\eq\tfrac{12}{35},
\leqno(3.5.1)$$
\par\nin where $\al'(h_1)\defs\alt(h_1)\mi\alt_f$, etc. Set $\om_0\eq \ddd x{\wedge}\ddd y{\wedge}\ddd z$, and
$$\be_k\defs\alt(\dd_t^k[h_1^k\om_0])\eq\tfrac{12+k}{35}\q(k\ins\N).$$
\par\nin In the notation of (\hl{1.1.5}{1.1.5}) we can verify that the degree $\be_r$ part
\htt{3.5.2}{}
$$(\dd_tt{-}\be_3)(\dd_tt{-}\be_2)(\dd_tt{-}\be_1)(\dd_t[h_1\om_0])_{\be_r}=\ct(\dd_t^r[h_1^r\om_0])_{\be_r}\in\G_f^{(\be_r)}
\leqno(3.5.2)$$
\par\nin coincides with $(\dd_t[h_2\om_0])_{\be_r}$ up to a nonzero constant multiple (where $\ct\ins\C^*$), using (\hl{3.5.1}{3.5.1}) and the inverse of (\hl{1.1.2}{1.1.2}) together with a generalization of (\hl{2.2.1}{2.2.1}), since
\htt{3.5.3}{}
$$\be_1\pl\al'(h_2){-}1\eq\be_5>\be_r.
\leqno(3.5.3)$$
\par\nin Note that $N\eq 0$, and the action of $\dd_tt$ is semisimple in (\hl{3.5.2}{3.5.2}). Here we consider the elements modulo $V^{>\be_r}\G_f$, when we calculate $(\dd_t^r[h_1^r\om_0])_{\be_r}$ using the inverse of (\hl{1.1.2}{1.1.2}) together with the first equality of (\hl{3.5.1}{3.5.1}). The divisions by $\dd_yf,\dd_zf$ needed in the definition of $\dd_t$ are then replaced by those by $\dd_yf_1,\dd_zf_1$, see also \cite{sem}.
\sk
By (\hl{3.5.2}{3.5.2}) we see that $[\om_0]_{\be_r}$ vanishes if we choose $c\ins\C^*$ appropriately. The degree $\be_r$ part $\Htppf{}^{(\be_r)}$ then vanishes, since the second smallest spectral number $\alt([x\om_0]$) is $\tfrac{12+5}{35}\sgt\be_r$. These show that Theorem~\hl{T2}{2}\,(ii) does not hold without the assumption on the image of $f_{\!\rho_f}^{\1 r}h$ with $r\eq 4$, and we can prove the assertion only for the minimal degree part of $f_{>1}$.
\par\htt{R3.5a}{}\msn
{\bf Remark 3.5a.} If we replace $x^7\pl y^5$ with $\tfrac{1}{7}\1x^7\pl\tfrac{1}{5}\1y^5$, one can verify that the constant $c$ in the definition of $f_{>1}$ is given by $\tfrac{6\cdot 8\cdot 3}{3\cdot 2\cdot 4}\eq 6$ as in \cite{Kat1} (replacing $x,y$ with $-x,-y$), see also \cite{sem}. Indeed, we have $\ct\eq\bl(\!\1{-}\tfrac{1}{35}\br)^{r-1}$ with $r\eq4$ in (\hl{3.5.2}{3.5.2}), since $c_{h,\be}\eq 1{-}\be$ in (\hl{2.2.1}{2.2.1}). So the left-hand side of (\hl{3.5.2}{3.5.2}) is given by $\tfrac{\!\!3\cdot 2}{35^3}(\dd_t[h_1\om_0])_{\be_r}$, and this contributes to $3\scd 2$ in the denominator, where 4 comes from $\al'(h_2)\mi 1\eq\tfrac{4}{35}$. Concerning the numerator $6\scd 8\scd 3$, this is obtained by derivations of monomials when we calculate $\dd_t^3$ using the inverse of (\hl{1.1.2}{1.1.2}) together with the first equality of (\hl{3.5.1}{3.5.1}), since $f_1$ is replaced as above so that $(f_1)_x\eq x^6$, $(f_1)_y\eq y^4$, see also \cite{sem}. (This is compatible with a computation using ``bernstein" in \cite[gmssing.lib]{Sing}, where we get a root $-\tfrac{51}{35}$.)
\par\htt{R3.5b}{}\msn
{\bf Remark 3.5b.} One may apply the above argument to calculate an example in \cite{Ca}, for instance, $f\eq x^7\pl y^6\pl x^5y^2\pl c\1x^3y^4$ with $c\eq\tfrac{4}{7}\scd \tfrac{2}{4}\eq\tfrac{2}{7}$. Indeed, we have $\alpha'(x^5y^2){-1}\eq\tfrac{2}{42}$, $\alpha'(x^3y^4){-1}\eq\tfrac{2r}{42}$, and $(x^5y^2)^r/x^7\eq x^3y^4$ with $r\eq 2$. The latter implies that $\tfrac{4}{7}$ comes from a calculation of $\dd_t$. We see that $\ct$ in the right-hand side of (\hl{3.5.2}{3.5.2}) is given by $-\tfrac{2}{42}$, which coincides up to sign with the coefficient in the left-hand side. Then $\tfrac{2}{4}$ comes from the ratio between $\alpha'(x^5y^2){-1}$ and $\alpha'(x^3y^4){-1}$ in view of a remark after (\hl{2.2.1}{2.2.1}) about $c_{h,\be}$.
\sk
One may also consider $f\eq x^7\pl y^6\pl x^5y^2\pl c\1x^3y^4\pl c'x^4y^4$. Here $c\eq\tfrac{5}{14}$, since we have to apply the above argument to the asymptotic expansion of $[x\om_0]$ (instead of $[\om_0]$). It is not easy to determine $c'$, since there are many nontrivial terms on the right-hand side of (\hl{2.2.1}{2.2.1}) and the calculation becomes much more complicated. As a conclusion it can be obtained up to sign by dividing 
$$\aligned&2^5\scd 19\scd 12\scd 5\scd 5/(7^3{\cdot}\1 6)\\
{}-{}&(2\pl 4\pl 6\pl 8)\scd 2^3\scd 4\scd 12\scd 5\scd 5\scd (5/14)/(7^2{\cdot}\1 6)\\
{}+{}&(8\scd 4\pl 8\scd 2\pl 6\scd 2)\scd 2\scd 4^2 \scd 5\scd 5\scd (5/14)^2/(7\scd 6)
\endaligned
$$
\par\nin by $8\scd 6\scd 4\scd 2\scd 10$, and the result is $5/16464\eq 5/(7^3{\cdot}\1 48)$, see also \cite{Ca}. Here $2\pl 4\pl 6\pl 8$ and $8\scd 4\pl 8\scd 2\pl 6\scd 2$ come from some combinatorial argument related to a generalization of (\hl{3.5.2}{3.5.2}), see also \cite{sem}. (This is compatible with a calculation using ``bernstein" in \cite[gmssing.lib]{Sing}, where we can get a root $-\tfrac{65}{42}$ if we put $c\eq(5/14)$, $c'\eq{-}(5/16464)$.)
\par\htt{3.6}{}\msn
{\bf 3.6.~Reduced geometric genus.} Let $(Z,0)$ be a hypersurface having an isolated singularity at $0$ with $d_Z\defs\dim Z\gess 2$. Let $\pi\cols\Zt\tos Z$ be a desingularization such that $E\defs\pi^{-1}(0)$ is a divisor with simple normal crossings and $\pi$ induces an isomorphism over $Z\stm\{0\}$.
\par\htt{P3.6}{}\msn
{\bf Proposition 3.6.} {\it We have the equality}
\htt{3.6.1}{}
$$\dim_{\C}H^{d_E}(E,\OO_E)=n_{f,1}.
\leqno(3.6.1)$$
\par\nin \msn
{\it Proof.} The mapping cone of the inclusion
\htt{3.6.2}{}
$$\Om_{\Zt}^{\ssb}(\log E)(-E)\into\Om_{\Zt}^{\ssb}(\log E)
\leqno(3.6.2)$$
\par\nin is essentially (or cohomologically) supported on $E$, and gives the Hodge complex of the {\it link} $L_{Z,0}\defs Z\,{\cap}\,S$ with $S$ a sufficiently small sphere in the ambient space $X$ forgetting the weight filtration, see also \cite{DuSa}. (Note that the complex $\Om_{\Zt}^{\ssb}(\log E)(-E)$ represents the 0-extension in $D^b_c(\Zt,\C)$, see for instance \cite[Prop.\,3.11\,(ii)]{mhm}.) Here the Hodge filtration is defined by the truncations $\sigma_{\ges p}$. By the strictness of the Hodge filtration (that is, the $E_1$-degeneration of the spectral sequence), this implies that
\htt{3.6.3}{}
$$\dim_{\C}H^{d_E}(E,\OO_E)=\dim_{\C}\Gr_F^0H^{d_E}(L_{Z,0},\C).
\leqno(3.6.3)$$
\par\nin \sk
On the other hand, it is known (see for instance \cite[(2.3.1)]{FPS}) that
\htt{3.6.4}{}
$$\dim_{\C}\Gr_F^0H^{d_E}(L_{Z,0},\C)=\dim_{\C}\Gr_F^1H^{d_Z}(F_{\!f},\C)_1,
\leqno(3.6.4)$$
\par\nin where $H^{d_Z}(F_{\!f},\C)_1$ denotes the unipotent monodromy part of the vanishing cohomology. (This is closely related to the Wang sequence, see \cite{Mi}.) Note that $\Gr_F^1H^{d_Z}(F_{\!f},\C)_1$ is contained in the kernel of $\Gr_FN$, since $\Gr_F^0H^{d_Z}(F_{\!f},\C)_1\eq0$, where we can also apply the symmetry of spectral numbers (\hl{1.2.3}{1.2.3}). By the definition of multiplicities of spectral numbers (\hl{1.2.1}{1.2.1}) we then get the equality (\hl{3.6.1}{3.6.1}) employing the symmetry of spectral numbers (\hl{1.2.3}{1.2.3}). This finishes the proof of Proposition~\hl{P3.6}{3.6}.
\par\htt{R3.6a}{}\msn
{\bf Remark 3.6a.} The left-hand side of (\hl{3.6.1}{3.6.1}) is called the {\it reduced genus\1} in \cite{BiSc}, \cite{MO}. By \cite{geo} there is the equality
\htt{3.6.5}{}
$$p_g(Z)=\msum_{\al\les 1}\,n_{f,\al},
\leqno(3.6.5)$$
\par\nin where the left-hand side is the {\it geometric genus\1} of $(Z,0)$.
\par\htt{R3.6b}{}\msn
{\bf Remark 3.6b.} We have a short exact sequence
\htt{3.6.6}{}
$$0\to\om_{\Zt}\to\om_{\Zt}(E)\to\om_E\to 0,
\leqno(3.6.6)$$
\par\nin applying the dual functor $\DD$ to the short exact sequence $0\tos\OO_{\Zt}(-E)\tos\OO_{\Zt}\tos\OO_E\tos 0$. By the Grauert-Riemenschneider vanishing theorem, it gives a short exact sequence
\htt{3.6.7}{}
$$0\to\pi_*\om_{\Zt}\to\pi_*\om_{\Zt}(E)\to\Gamma(E,\om_E)\to0.
\leqno(3.6.7)$$
\par\nin As is well known, there is a canonical inclusion $\pi_*\om_{\Zt}(E)\into\om_Z$, which induces
\htt{3.6.8}{}
$$\Gamma(E,\om_E)\into(\om_Z/\pi_*\om_{\Zt})_0.
\leqno(3.6.8)$$
\par\nin Here the left-hand side is the dual of $H^{d_E}(E,\OO_E)$ by Grothendieck duality, and the dimension of the right-hand side is the geometric genus.
\par\htt{3.7}{}\msn
{\bf 3.7.~Non-rational Du Bois singularity case.} It is known that the multiplicity of the minimal spectral number is 1, see for instance \cite[\S 4.11]{DS1}. Hence we have $n_{f,1}\eq1$ if $Z$ has a non-rational Du Bois singularity at 0 (that is, $\alt_f\eq 1$, see \cite[Thm.\,0.4]{rat}, \cite[Thm.\,0.5]{mos}. In this case we have a positive answer to the conjecture \cite[Conj.\,1.7]{BiSc} using a property of (B)-lattices explained in \cite[\S 3.1]{bl}. This is compatible with \cite[Cor.\,1.2]{MO}. (It does not seem easy to verify \cite[Thm.\,1.1]{MO} employing our method.) These justify the definition (\hl{3}{3}), which considers the subspaces modulo the image of $N$ when $N\nes0$. (This comes from the {\it Verdier-type extension theorem\1} as is explained in the proof of Theorem~\hl{T1}{1}.)
\sk
In the case $f\eq x_1\cdots x_n\pl\msum_{i=1}^n\,x^{b_i}$ with $\msum_{i=1}^n\,1/b_i\slt 1$, for instance, we have $\alt_f\eq 1$ (see for instance \cite{exp}) and $N^{n-2}\nes 0$ on $\G_f^{(1)}$ as a corollary of the descent theorem of nearby cycle formula \cite[Cor.\,3]{des}. Hence the multiplicity of the root $-1$ of $b_f(s)/(s{+}1)$ is $n{-}1$, and
\htt{3.7.1}{}
$$\dim_{\C}\Htppf{}^{(1)}=n\mi1,
\leqno(3.7.1)$$
\par\nin since $\Hppf$ is a (B)-lattice in the sense of \cite[\S 3.1]{bl}. (Note that $N$ is given by the action of $\dd_tt{-}\al\eq{-}s{-}\al$ on $\Gr_V^{\al}$.) As for the integral spectral numbers, we have $n_{f,j}\eq 1$ if $j\eq 1$ or $n{-}1$, and $n_{f,j}\gess 1$ if $j\in[2,n{-}2]$. We can determine the $n_{f,\al}$ using a combinatorial formula in \cite{St}, see also \cite{JKSY1}.
\bs\bs
\vbox{\centerline{\bf 4. Simpler proof of Theorem~\hl{T1}{1}}
\bsn
In this section we explain a simple proof of a generalization of the last assertion in \cite[Cor.\,3.0.1]{Bit} to the non-integral exponent case implying a simpler proof of Theorem~\hl{T1}{1}.}
\par\htt{4.1}{}\msn
{\bf 4.1.~Another approach.} Let $f$ be a holomorphic function on a complex manifold $X$ of dimension $n\gess 2$ having a unique singular point at $0\ins X$. Set $Z\defs f^{-1}(0)\sst X$. Let $r_{\!f}$ be the number of local irreducible components at 0. For $\be\ins\Q\,{\cap}\,[0,1)$, let $P$ be the {\it pole order\1} filtration on $\OO_X(*Z)f^{\be}\defs\OO_X[f^{-1}]f^{\be}$ such that
$$P_k\bl(\OO_X(*Z)f^{\be}\br)\defs\OO_Xf^{\be-1-k}\q\h{if}\,\,\,k\ges0,$$
\par\nin and 0 otherwise. Let $\PD$ be the filtration on $\OO_X(*Z)f^{\be}$ defined by
$$\PD_k\bl(\OO_X(*Z)f^{\be}\br)\defs\D_Xf^{\be-1-k}\q\h{if}\,\,\,k\ges0,$$
\par\nin and 0 otherwise. Set
$$H^{(\be)}\defs\Hc^0\DR_X\bl(\OO_X(*Z)f^{\be}\br){}_0.$$
\par\nin This is isomorphic to $\Hc^0i_0^*\RR j_*L^{(\be)}[n]$ with $L^{(\be)}$ the pull-back of a local system of rank 1 on $\C^*$ with monodromy $e^{-2\pi i\be}$ by $f$ and
$i_0\cols\{0\}\,{\into}\,X$, $j\cols X\stm Z\,{\into}\,X$ canonical inclusions. 
It has a canonical mixed complex Hodge structure. (Sometimes $f^{\be}$ is omitted in the case $\be\eq0$.)
\sk
The filtrations $P,\PD$ can induce (quotient) filtrations on $H^{(\be)}$, which are denoted also by $P,\PD$ respectively. Recall that the de Rham functor $\DR_X$ is shifted by $n$ to the left for the Riemann-Hilbert correspondence as usual. Choosing local coordinates $x_1,\dots,x_n$, we have
$$\aligned H^{(\be)}&\eq\OO_{X,0}(*Z)f^{\be}/\bl(\msum_{i=1}^n\,\dd_{x_i}\bl(\OO_{X,0}(*Z)f^{\be}\br)\br),\\ \D_X&\eq\OO_X\pl\msum_{i=1}^n\,\dd_{x_i}\D_X.\endaligned$$
\par\nin (Here it may be more natural to use the de Rham functor for right $\D$-modules.) Using the functional equation in the definition of $b_f(s)$ (or by the holonomicity), we have
\htt{4.1.1}{}
$$\OO_X(*Z)f^{\be}=\D_Xf^{\be-1-k_0}\q\h{for {\it any}}\,\,k_0\,{\gg}\,0.
\leqno(4.1.1)$$
\par\nin These immediately imply the following.
\par\htt{L4.1}{}\msn
{\bf Lemma~4.1.} {\it For $\be\ins\Q\,{\cap}\,[0,1)$, we have the equalities}
$$P_kH^{(\be)}\eq\PD_kH^{(\be)}\,\,\,\,\h{\it in}\,\,\,H^{(\be)}\q(\forall\,k\ins\Z).$$
\par\nin \ms
In the case $\be\eq0$, the algebraic local cohomology sheaf
$$\M_Z^{\vee}\defs\Hc^1_{[Z]}(\OO_X)\eq\OO_X(*Z)/\OO_X$$
\par\nin corresponds to the dual of the constant sheaf $\C_Z[n{-}1]$ up to a Tate twist by the de Rham functor $\DR_X$, and it contains a regular holonomic $\D_X$-module $\M_Z^{\rm IC}$ corresponding to the intersection complex ${\rm IC}_Z\C$. Set
$$\M''\defs\M_Z^{\vee}/\M_Z^{\rm IC}\eq\OO_X(*Z)/\M',$$
\par\nin where $\M'$ is defined by the last isomorphism, that is, the kernel of the canonical surjection, which is an extension of $\M_Z^{\rm IC}$ by $\OO_X$. We see that $\M''$ is supported at $0$, and is isomorphic to a finite direct sum of copies of a simple $\D_{X,0}$-module $\B_0\defs\C[\dd_{x_1},\dots,\dd_{x_n}]$. We have
\htt{4.1.2}{}
$$\Hc^0\DR_X(\M')=0,
\leqno(4.1.2)$$
\par\nin since
$$\Hc^0\DR_X(\OO_X)\eq\Hc^0\DR_X(\M_Z^{\rm IC})\,\bl({=}\,\Hc^0{\rm IC}_X\C\br)\eq0.$$
\par\nin \sk
In the case $\be\ins(0,1)$, we have a simple $\D_X$-submodule $\M'\sst\OO_X(*Z)f^{\be}$ whose quotient $\M''\defs\OO_X(*Z)f^{\be}/\M'$ is supported on $0$, since $\OO_X(*Z)f^{\be}$ contains no nontrivial coherent $\D_X$-submodule supported at 0. The vanishing (\hl{4.1.2}{4.1.2}) also holds, since $\M'$ corresponds to the intersection complex ${\rm IC}_XL^{(\be)}$.
\sk
For any $\be\ins[0,1)$, consider the induced and quotient filtrations $\PD$ on $\M',\M''$. These are filtrations by $\D_X$-{\it submodules.} Apply the snake lemma to the commutative diagram below after replacing the first term of the upper row with an appropriate quotient so that the row becomes a short exact sequence:
\htt{4.1.3}{}
$$\begin{array}{ccccccccccc}
\Hc^0\DR_X(\PD_k\M')&\!\!\!\to&\!\!\!\Hc^0\DR_X(\PD_k\OO_X(*Z)f^{\be})&\!\!\!\to&\!\!\!\Hc^0\DR_X(\PD_k\M'')&\!\!\!\to&\!\!\!0\\
\downarrow&&\,\,\downarrow\raise1pt\h{$\!\scriptstyle\beta_k$}&\raise5mm\h{}\raise-3mm\h{}&\,\,\downarrow\raise1pt\h{$\!\scriptstyle\gamma_k$}\\
0&\!\!\!\to&\!\!\!\Hc^0\DR_X(\OO_X(*Z)f^{\be})&\!\!\!\to&\!\!\!\Hc^0\DR_X(\M'')&\!\!\!\to&\!\!\!0\end{array}
\leqno(4.1.3)$$
\par\nin We then get the isomorphism ${\rm Coker}\,\beta_k\simto{\rm Coker}\,\gamma_k$, and hence ${\rm Im}\,\beta_k\simto{\rm Im}\,\gamma_k$, where ${\rm Im}\,\beta_k\eq P^{\D}_kH^{(\be)}$ (since the latter is essentially a quotient filtration). So the following is proved.
\par\htt{C4.1}{}\msn
{\bf Corollary~4.1.} {\it For $\be\ins\Q\,{\cap}\,[0,1)$, we have the equalities
$$\ell_{\D_{X,0}}(\D_{X,0}f^{\be-1-k})=\dim_{\C}P_kH^{(\be)}+r_{\!f}\delta_{\be,0}+1\q(\forall\,k\ins\Z_{\ges0}),$$
\par\nin where the left-hand side is the length of $\D_{X,0}f^{\be-1-k}$ as a $\D_{X,0}$-module.}
\ms
In the case $\be\eq0$, $n\gess3$ this has been obtained by T.\,Bitoun (see the last assertion in \cite[Cor.\,3.0.1]{Bit}) employing a much more complicated argument. For the proof of the above assertion, we do not need a formula of Vilonen determining the intersection complex $\D$-module. We can see that Corollary~\hl{C4.1}{4.1} implies a simpler proof of Theorem~\hl{T1}{1}, see \hl{4.2}{4.2} below. It does not seem however easy to apply Corollary~\hl{C4.1}{4.1} to explicit computations for the case $\be\nes0$ (without hiring an argument similar to Remark~\hl{R4.1c}{4.1c} below).
\sk
If $\be\ins(0,1)$, the above proof of Corollary~\hl{C4.1}{4.1} holds for any $k\ins\Z$. Indeed, for any $\be\ins[0,1)$, both sides of the equality coincide with 1 if $k\ins\Z_{<0}$. One can easily see that the left-hand side is 1 for any $\be\ins[0,1)$, $k\ins\Z_{<0}$, since $\D_Xf^{\al}\eq\D_Xf^{\al+1}$ for $b_f(\al)\nes0$ and the image of $f^{\be}$ in $\M''$ is annihilated by $f^j$ for $j\,{\gg}\,0$.
\par\htt{R4.1a}{}\msn
{\bf Remark~4.1a.} Assume $\be\eq0$. Let $\psi_{f,1}\C_X[n']$, $\varphi_{f,1}\C_X[n']$ be the unipotent monodromy part of the nearby and vanishing cycle complexes, where $n'\defs n{-}1$. Let $H$ be the vanishing cohomology group of $f$ at 0 with $H_1,H_{\ne1}\sst H$ its unipotent and non-unipotent monodromy parts which are identified with $\varphi_{f,1}\C_X[n']$ and $\varphi_{f,\ne1}\C_X[n']$ respectively. Set $N\defs\log T_u$ with $T\eq T_sT_u$ the Jordan decomposition of the monodromy $T$. We have the canonical isomorphism (using for instance \cite[(3.1.3)]{JKSY3})
$${\rm Coker}\bl(N\cols\psi_{f,1}\C_X[n']\tos\psi_{f,1}\C_X(-1)[n']\br)=\DR_X\bl(\Hc^1_{[Z]}(\OO_X)\br),$$
\par\nin where ${\rm Ker}\,N$ is defined in the abelian subcategory of $D^b_c(X,\C)$ associated with the middle perversity in \cite{BBD}. The right-hand side of the isomorphism is identified with
$$\RR\Gamma_{Z}\C_X[n{+}1]\eq(i_Z)_*i_Z^!\C_X[n'{+}2],$$
\par\nin (where $i_Z\cols Z\,{\into}\,X$ denotes the inclusion), and also with the dual of $\C_Z[n']$ up to a Tate twist. (Note that $H^{(0)}$ is denoted by $H'$ in \cite{Bit}.)
Using (\hl{4.1.2}{4.1.2}) together with the $N$-{\it primitive decompositions\1} of $\Gr^W_{\ssb}\psi_{f,1}\C_X[n']$, $\Gr^W_{\ssb}\varphi_{f,1}\C_X[n']$ (where their primitive parts coincide for weights at least $n$, see for instance \cite[(2.2.5)]{KLS}), we can then get the isomorphism of mixed Hodge structures
\htt{4.1.4}{}
$$H^{(0)}\eq H_1/NH_1.
\leqno(4.1.4)$$
\par\nin So Corollary~\hl{C4.1}{4.1} is actually {\it quite similar\1} to Theorem~\hl{T1}{1}.
\sk
The assertion (\hl{4.1.4}{4.1.4}) is closely related to the {\it Wang sequence\1} in \cite{Mi}. Note however that the {\it topological variation\1} $T{-}\1{\rm id}$ is {\it not\1} a morphism of mixed Hodge structures unless the monodromy $T$ is {\it semisimple\1} or $(T{-}\1{\rm id})^2\eq0$. (One {\it cannot\1} add morphisms of mixed Hodge structures with {\it different types.}) One has to prove that its image {\it coincides with\1} $NH_1\oplus H_{\ne1}$ in order to show that its cokernel has a canonical mixed Hodge structure.
\sk
The isomorphism (\hl{4.1.4}{4.1.4}) is also compatible with the induced Hodge filtration $F$ because of (\hl{4.1.2}{4.1.2}). (This is closely related to {\it primitive spectral numbers\1} which is defined by using the spectral pairs.) This assertion does {\it not\1} follow from the {\it general\1} theory of pull-back functors of Hodge modules, since the {\it naively\1} induced Hodge filtration on $H^{(0)}$ is {\it not\1} a {\it good\1} one in general. Here we need in an essential way the assertion that the Hodge filtration on $\M''$ coincides with the quotient filtration of the Hodge filtration $F$ (that is, the {\it strictness\1} of the Hodge filtration). Note also that $F\sst P$ by definition.
\par\htt{R4.1b}{}\msn
{\bf Remark~4.1b.} Assume $\be\ins(0,1)$. Let $\psi_{f,\la}\C_X[n']$, $\varphi_{f,\la}\C_X[n']$ be the $\la$-eigenpart of the nearby and vanishing cycle complexes with $\la\defs e^{2\pi i\be}$. The latter is identified with the (generalized) $\la$-eigenspace of the vanishing cohomology $H_{\la}\sst H$. Since $\la\nes 1$, we have the canonical isomorphism
$$\psi_{f,1}L^{(\be)}[n']=\psi_{f,\la}\C_X[n']=\varphi_{f,\la}\C_X[n']=H_{\la}.$$
\par\nin We apply the functorial isomorphism
$$i_Z^*=C({\rm can}\cols\psi_{f,1}\tos\varphi_{f,1})[-1]$$
\par\nin to $\RR j_*L^{(\be)}[n]$ with $j\cols X\stm Z\,{\into}\,X$ the inclusion. Here the canonical morphism
$${\rm can}\cols\psi_{f,1}L^{(\be)}[n']\tos\varphi_{f,1}\RR j_*L^{(\be)}[n']$$
\par\nin is identified with $N\cols\psi_{f,1}L^{(\be)}[n']\tos\psi_{f,1}L^{(\be)}(-1)[n']$, since $N\eq{\rm can}\ssc{\rm Var}$ and
$${\rm Var}\cols\varphi_{f,1}\RR j_*L^{(\be)}[n']\to\psi_{f,1}L^{(\be)}(-1)[n']$$
\par\nin is an isomorphism. (Indeed, the morphism Var is defined by the action of $t$ up to a sign for the corresponding $\D$-module using the $V$-filtration and the graph embedding.) We thus get the isomorphism
\htt{4.1.5}{}
$$H^{(\be)}=H^0i_0^*\RR j_*L^{(\be)}[n]={\rm Coker}\bl(N\cols H_{\la}\tos H_{\la}(-1)\br).
\leqno(4.1.5)$$
\par\nin So Corollary~\hl{C4.1}{4.1} is quite similar to Theorem~\hl{T1}{1} also in the case $\be\ins(0,1)$.
\par\htt{R4.1c}{}\msn
{\bf Remark~4.1c.} For $\al\ins\Q$, set
$$\aligned\eta_{\al,i}&\defs f\dd_i\pl\al f_i\in\D_{X,0}\q\h{with}\q\dd_i\defs\dd_{x_i},\,\,\,f_i\defs\dd_if\,\,\,\,(i\ins[1,n]),\\ E^{\lng\al\rng}&\defs\rlap{\raise-10.5pt\h{$\,\,\,\scriptstyle k$}}\indlim\,E_{\al-k}\q\h{with}\q E_{\al-k}\defs\OO_{X,0}/\msum_{i=1}^n\,\eta_{\al-k,i}\1\OO_{X,0}\,\,\,(k\ins\N),\endaligned$$
\par\nin where the transition morphism $E_{\al-k}\tos E_{\al-k-1}$ is induced by multiplication by $f$ using the relation $f\ssc\eta_{\al-k,i}\eq\eta_{\al-k-1,i}\ssc f$. By definition we have the isomorphism
\htt{4.1.6}{}
$$E^{\lng\al\rng}=\Hc^0\DR_X\bl(\OO_X(*Z)f^{\al}\br){}_0.
\leqno(4.1.6)$$
\par\nin \sk
On the other hand, there is the isomorphism
$$\Hppf[t^{-1}]=\rlap{\raise-10.5pt\h{$\,\,\,\scriptstyle k$}}\indlim\,\E_k\q\h{with}\q \E_k\defs\Hppf\,\,\,(k\ins\N),$$
\par\nin where the transition morphism is given by the action of $t$ (induced by multiplication by $f$). We see that
$$\eta_i\ssc f_j\mi\eta_j\ssc f_i\eq f\ssc(f_i\ssc\dd_j\mi f_j\ssc\dd_i)\eq(f_i\ssc\dd_j\mi f_j\ssc\dd_i)\ssc f,$$
\par\nin hence there is a surjective morphism
\htt{4.1.7}{}
$$\Hppf[t^{-1}]\onto E^{\lng\al\rng},
\leqno(4.1.7)$$
\par\nin by the definition of the Brieskorn lattice, see (\hl{1.1.1}{1.1.1}). This induces the surjection
\htt{4.1.8}{}
$$\Hppf[t^{-1}]/(\dd_tt\pl\al)\Hppf[t^{-1}]\onto E^{\lng\al\rng},
\leqno(4.1.8)$$
\par\nin by the definition of the Gauss-Manin connection (\hl{1.1.2}{1.1.2}) using the inductive limit expression. This surjection is an isomorphism by (\hl{4.1.4}{4.1.4}--\hl{4.1.5}{4}). This argument is led to the formulation in \hl{4.2}{4.2} below.
\par\htt{4.2}{}\msn
{\bf 4.2.~Simpler proof of Theorem~\hl{T1}{1}.} We can show that Corollary~\hl{C4.1}{4.1} implies a simple proof of Theorem~\hl{T1}{1} employing the commutative diagram
\htt{4.2.1}{}
$$\begin{array}{cccccccccccccccc}&&\!\!\Hc^0\DR_X(\D_X[s]f^s)&\!\!\onto&\!\!\Hc^0\DR_X(\D_Xf^{-\al})\\ &&\downarrow&&\downarrow\\ \!\!\Hc^0\DR_X(\D_X[s]f^s)&\!\!\buildrel{s+\al'}\over{\longrightarrow}&\!\!\Hc^0\DR_X(\D_X[s]f^s)&\!\!\onto&\!\!\Hc^0\DR_X(\D_Xf^{-\al'})\end{array}
\leqno(4.2.1)$$
\par\nin (using the right-exactness of $\Hc^0\DR_X$) together with (\hl{4.1.1}{4.1.1}) and also the isomorphism
\htt{4.2.2}{}
$$\Hc^0\DR_X(\D_X[s]f^s)\eq\Htppf\,\,\bl(=\widehat{\mopl}_{\al\in\Q}\,\Htppf{}^{(\al)}\br),
\leqno(4.2.2)$$
\par\nin see \cite{Ma} (and (\hl{1.1.4}{1.1.4}) for the last equality). Here the middle vertical morphism is induced by $t^k$  (sending $f^s$ to $f^{s+k}$ and $s$ to $s{+}k$, see \cite{Ka}) with $k\defs\al'\mi\al\ins\Z$. We assume that $\al'$ is sufficiently large so that there is a short exact sequence
\htt{4.2.3}{}
$$0\to\D_X[s]f^s\buildrel{s{+}\al'}\over\longrightarrow\D_X[s]f^s\to\D_Xf^{-\al'}\to0,
\leqno(4.2.3)$$
\par\nin see \cite[Prop.\,6.2]{Ka} (and also \cite[Prop.\,A]{BaSa}). Recall that $s$ is identified with $-\dd_tt$.
\par\htt{R4.2a}{}\msn
{\bf Remark~4.2a.} The assertion of \cite[Cor.\,3.0.4]{Bit} should be divided into two parts: The first one is about GAGA, and is quite standard. The second part should contain the {\it comparison\1} between $H^n_{\rm dR}(B\stm Y^{\rm an})$ and $H^n_{\rm dR}\bl(B,\OO_B(*Y^{\rm an}_B)\br)$, where $B$ is a sufficiently small  analytic open ball with $Y^{\rm an}_B\defs Y^{\rm an}\cap B$. This is {\it purely analytic,} although it is essentially the Grothendieck theorem. It remains to show the independence of the radius of $B$, but this is also entirely analytic. (It can be proved by using only embedded resolutions of singularities for the closures of strata.)
\par\htt{R4.2b}{}\msn
{\bf Remark~4.2b.} In \cite[Cor.\,3.0.4]{Bit}, it is not very clear whether the assertion is proved in the general situation, since it does not seem to be shown that for any $f\ins\OO_{Y,0}$ with $Y$ a smooth complex algebraic variety, there is a local coordinate system $(y_1,\dots,y_n)$ of $Y$ around 0 such that $f$ is a polynomial of $y_1,\dots,y_n$. Actually one can apply GAGA to the {\it general setting\1} without relying on the above assertion.
\sk
The latter assertion does not necessarily hold in general. Consider for instance the elliptic curve case with $n\eq1$. We may assume that $f$ and $g\defs y_1$ are extended to finite morphisms $\ft\,{:}\,\Yt\tos\PP^1$, $\gt\,{:}\,\Yt\tos\PP^1$ taking the normalization $\Yt$ of a compactification of $Y$ (and applying Riemann's extension theorem). We have
\htt{4.2.4}{}
$${\rm div}\,\ft\eq\msum_{P\in\Yt}\,m_P\1[P],\q{\rm div}\,\gt\eq\msum_{P\in\Yt}\,m'_P\1[P],
\leqno(4.2.4)$$
with $m_P,m'_P\ins\Z$ and $m'_0\eq1$. Assume the following relation holds on a non-empty analytic open subset of $\Yt$\,:
\htt{4.2.5}{}
$$\ft\eq\msum_{i=k}^{k'}\,c_i\1\gt^i,
\leqno(4.2.5)$$
where $c_i\ins\C$, $c_kc_{k'}\nes0$. Then it holds on the complement $U$ of the union of the poles of $\ft$ and $\gt$. We then get that
\htt{4.2.6}{}
$$m_P=k\1m'_P,\,\,\,\h{hence}\,\,\,\,m_P\ins m_0\Z\q(\forall\,P\ins\gt^{-1}(0)),
\leqno(4.2.6)$$
since $\gt(0)\eq0$ and $m'_0\eq1$. Moreover the poles of $\ft$ coincide with those of $\gt$. (Indeed, if $|\gt(y)|$ is sufficiently large, we have $\msum_{i=k}^{k'-1}\,|c_i\1\gt(y)^i|\less\tfrac{1}{2}|c_{k'}\gt(y)^{k'}|$, hence $|\ft(y)|\gess\tfrac{1}{2}|c_{k'}\gt(y)^{k'}|$. On the other hand, $|\ft(y)|$ is bounded by $C|\gt(y)|^{k'}$ or $C|\gt(y)|^{k}$ for $y\ins U$ depending on whether $|\gt(y)|\gess1$ or $|\gt(y)|\less1$, where $C\defs\msum_{i=k}^{k'}\,|c_i|$.) There are many rational functions $f$ on elliptic curves such that the last inclusion of (\hl{4.2.6}{4.2.6}) cannot be satisfied. One may consider for instance the case $m_0\sgt1$ and $m_P\eq1$ for any $P\ins\ft^{-1}(0)\stm\{0\}$. (On an elliptic curve, there is a rational function $\ft$ satisfying the first equality of (\hl{4.2.4}{4.2.4}) if $\msum_P\,m_PP\eq0$, choosing the point of the complement of $U$ as the origin.) This implies a contradiction, since $\gt^{-1}(0)\nes\{0\}$. (It does not seem trivial to generalize this argument to the higher-dimensional case, since the {\it poles of coordinates\1} are not easy to control.)

\sk
{\smaller\smaller RIMS Kyoto University, Kyoto 606-8502 Japan}

\begin{thebibliography}{JKSY\,22b}
\bibitem[BaSa\,22]{BaSa} Bath, D., Saito, M., Twisted logarithmic complexes of positively weighted homogeneous divisors (arxiv:2203.11716).
\bibitem[BBD\,82]{BBD} Beilinson, A., Bernstein, J., Deligne, P., Faisceaux pervers, Ast\'erisque 100, Soc.\ Math.\ France, Paris, 1982.
\bibitem[Bi\,23]{Bit} Bitoun, T., On the $D$-module of an isolated singularity (arxiv:2307.00120).
\bibitem[BiSc\,18]{BiSc} Bitoun, T., Schedler, T., On $\D$-modules related to the $b$-function and Hamiltonian flow, Compos.\ Math.\ 154 (2018), 2426--2440.
\bibitem[Br\,70]{Br} Brieskorn, E., Die Monodromie der isolierten Singularit\"aten von Hyperfl\"achen, Manuscripta Math., 2 (1970), 103--161.
\bibitem[Ca\,87]{Ca} Cassou-Nogu\`es, P., Etude du comportement du polyn\^ome de Bernstein lors d'une d\'eformation \`a $\mu$ constant de $X^a+Y^b$ avec $(a,b)=1$,
Compos.\ Math. 63 (1987), 291--313.
\bibitem[DGPS\,20]{Sing} Decker, W., Greuel, G.-M., Pfister, G., Sch\"onemann, H., {\sc Singular} 4.2.0 --- A computer algebra system for polynomial computations, available at http://www.singular.uni-kl.de (2020).
\bibitem[DiSa\,12]{DS1} Dimca, A., Saito, M., Koszul complexes and spectra of projective hypersurfaces with isolated singularities (arxiv:1212.1081).
\bibitem[DiSa\,14]{DS2} Dimca, A., Saito, M., Some remarks on limit mixed Hodge structures and spectrum, An.\ \c{S}t.\ Univ.\ Ovidius Constan\c{t}a Ser.\ Mat.\ 22 (2014), 69--78.
\bibitem[DuSa\,90]{DuSa} Durfee, A.\,H., Saito, M., Mixed Hodge structures on the intersection cohomology of links, Compos.\ Math.\ 76 (1990), 49--67.
\bibitem[FPS\,21]{FPS} Fern\'andez de Bobadilla, J., Pallar\'es, I., Saito, M., Hodge modules and cobordism classes (arxiv:2103.04836), to appear in JEMS.
\bibitem[GLS\,07]{GLS} Greuel, G.-M., Lossen, C., Shustin, E., Introduction to Singularities and Deformations, Springer, Berlin, 2007.
\bibitem[Gr\,61]{Gro} Grothendieck, A., El\'ement de g\'eom\'etrie alg\'ebrique III-1, Publ. Math. IHES 11 (1961), 5--167.
\bibitem[Ha\,75]{Ha} Hartshorne, R., On the de~Rham cohomology of algebraic varieties, Publ.\ Math.\ IHES 45 (1975), 5--99.
\bibitem[JKSY\,19]{JKSY1} Jung, S.-J., Kim, I.-K., Saito, M., Yoon, Y., Spectrum of non-degenerate functions with simplicial Newton polyhedra (arxiv:1911.09465).
\bibitem[JKSY\,22a]{JKSY2} Jung, S.-J., Kim, I.-K., Saito, M., Yoon, Y., Hodge ideals and spectrum of isolated hypersurface singularities, Ann.\ Inst.\ Fourier 72 (2022), 465--510 (Theorem and formula numbers are changed by the publisher, see arxiv:1904.02453 for the original numbers).
\bibitem[JKSY\,22b]{JKSY3} Jung, S.-J., Kim, I.-K., Saito, M., Yoon, Y., Higher Du Bois singularities of hypersurfaces, Proc.\ London Math.\ Soc.\ 125 (2022), 543--567.
\bibitem[Kas\,76]{Ka} Kashiwara, M., $B$-functions and holonomic systems, Inv.\ Math.\ 38 (1976/77), 33--53.
\bibitem[Kat\,81]{Kat1} Kato, M., The $b$-function of $\mu$-constant deformation of $x^7+y^5$, Bull.\ College Sci., Univ.\ Ryukyus 32 (1981), 5--10.
\bibitem[Kat\,82]{Kat2} Kato, M., The $b$-function of $\mu$-constant deformation of $x^9+y^4$, Bull.\ College Sci., Univ.\ Ryukyus 33 (1982), 5--8.
\bibitem[KLS\,22]{KLS} Kerr, M., Laza, R., Saito, M., Deformation of rational singularities and Hodge structure, Algebraic Geometry 9 (2022), 476--501.
\bibitem[Ma\,75]{Ma} Malgrange, B., Le polyn\^ome de Bernstein d'une singularit\'e isol\'ee, Lect.\ Notes in Math.\ 459, Springer, Berlin, 1975, pp. 98--119.
\bibitem[Mi\,68]{Mi} Milnor, J., Singular Points of Complex Hypersurfaces, Annals of Math.\ Studies No. 61, Princeton University Press, 1968.
\bibitem[MuOl\,22]{MO} Musta\c{t}\u{a}, M., Olano, S., On a conjecture of Bitoun and Schedler (arxiv:2207.02047).
\bibitem[SaK\,71]{SaK} Saito, K., Quasihomogene isolierte Singularit\"aten von Hyperfl\"achen, Inv.\ Math.\ 14 (1971), 123--142.
\bibitem[Sa\,83a]{geo} Saito, M., On the exponents and the geometric genus of an isolated hypersurface singularity, Proc. Sympos. Pure Math., 40, AMS, Providence, RI, 1983, pp.~465--472.
\bibitem[Sa\,83b]{sup} Supplement to: ``Gauss-Manin system and mixed Hodge structure", Ast\'erisque 101-102, Soc.\ Math,\ France, 1983, 320--331.
\bibitem[Sa\,84]{fil} Saito, M., Hodge filtrations on Gauss-Manin systems I, J.\ Fac.\ Sci.\ Univ.\ Tokyo Sect.\ IA Math.\ 30 (1984), 489--498.
\bibitem[Sa\,88a]{exp} Saito, M., Exponents and Newton polyhedra of isolated hypersurface singularities, Math.\ Ann.\ 281 (1988), 411--417.
\bibitem[Sa\,88b]{mhp} Saito, M., Modules de Hodge polarisables, Publ.\ RIMS, Kyoto Univ.\ 24 (1988), 849--995.
\bibitem[Sa\,89]{bl} Saito, M., On the structure of Brieskorn lattice, Ann.\ Inst.\ Fourier 39 (1989), 27--72.
\bibitem[Sa\,90]{mhm} Saito, M., Mixed Hodge modules, Publ.\ RIMS, Kyoto Univ.\ 26 (1990), 221--333.
\bibitem[Sa\,93]{rat} Saito, M., On $b$-function, spectrum and rational singularity, Math.\ Ann.\ 295 (1993), 51--74.
\bibitem[Sa\,94]{mic} Saito, M., On microlocal $b$-function, Bull.\ Soc.\ Math.\ France 122 (1994), 163--184.
\bibitem[Sa\,09]{mos} Saito, M., On the Hodge filtration of Hodge modules, Moscow Math.\ J.\ 9 (2009), 161--191.
\bibitem[Sa\,20]{des} Saito, M., Descent of nearby cycle formula for Newton non-degenerate functions (preprint, arxiv:2004.12367).
\bibitem[Sa\,21]{rp} Saito, M., $\D$-modules generated by rational powers of holomorphic functions, Publ.\ RIMS, Kyoto Univ.\ 57 (2021), 867--891.
\bibitem[Sa\,22a]{rh} Saito, M., Notes on regular holonomic $\D$-modules for algebraic geometers (arxiv:2201.01507).
\bibitem[Sa\,22b]{sem} Saito, M., Bernstein-Sato polynomials of semi-weighted-homogeneous polynomials of nearly Brieskorn-Pham type (arXiv:2210.01028).
\bibitem[ScSt\,85]{SS} Scherk, J., Steenbrink, J.\,H.\,M., On the mixed Hodge structure on the cohomology of the Milnor fibre, Math.\ Ann.\ 271 (1985), 641--665.
\bibitem[St\,77]{St} Steenbrink, J.\,H.\,M., Mixed Hodge structure on the vanishing cohomology, in Real and complex singularities, Sijthoff and Noordhoff, Alphen aan den Rijn, 1977, pp. 525--563.
\bibitem[Va\,81]{Va1} Varchenko, A.\,N., Asymptotic mixed Hodge structure in vanishing cohomologies, Izv.\ Akad.\ Nauk SSSR Ser.\ Mat.\ 45 (1981), 540--591.
\bibitem[Va\,82a]{Va2} Varchenko, A.\,N., The complex singularity index does not change along the stratum $\mu=$ const, Funk.\ Anal.\ Pri. 16 (1982), 1--12.
\bibitem[Va\,82b]{Va3} Varchenko, A.N., A lower bound for the codimension of the stratum $\mu=$ constant in terms of the mixed Hodge structure, Vestnik Moskow Univ.\ Ser.\ I Mat.\ Mekh.\ 37 (1982), 28--31. 
\end{thebibliography}
\end{document}